\documentclass[a4paper,twoside,reqno,11pt,tbtags]{article}

\RequirePackage{amssymb,amsmath,amsthm,mathrsfs}

\numberwithin{equation}{section}

\theoremstyle{plain}
\newtheorem{prop}{Proposition}[section]
\newtheorem{thm}[prop]{Theorem}
\newtheorem{lem}[prop]{Lemma}
\newtheorem{cor}[prop]{Corollary}
\theoremstyle{definition}
\newtheorem*{defi}{Definition}

\newtheorem{rem}[prop]{Remark}

\theoremstyle{remark}
\newtheorem{example}{Example}
\newtheorem*{msetting1}{$\Diamond$ Model Setting 1}
\newtheorem*{msetting1p}{$\Diamond$ Model Setting 1$'$}
\newtheorem*{msetting2}{$\Diamond$ Model Setting 2}

\DeclareMathOperator{\cov}{cov}

\DeclareMathOperator{\poisson}{Po}

\DeclareMathOperator{\interior}{int}

\DeclareMathOperator*{\esssup}{ess \, sup}

\DeclareMathOperator*{\sign}{sign}

\newcommand{\eps}{\varepsilon}

\newcommand{\bbeta}{\breve{\beta}}

\newcommand{\bbetaind}{\bbeta^{(\text{\rm ind})}}
\newcommand{\bbetasup}{\bbeta^{(\text{\rm sup})}}
\newcommand{\bgamma}{\breve{\gamma}}

\newcommand{\bphi}{\breve{\phi}}

\newcommand{\fbar}{\bar{f}}

\newcommand{\rbar}{\bar{r}}

\newcommand{\inlawto}{\, \stackrel{\text{\raisebox{-1pt}{$\mathscr{D}$}}}{\longrightarrow} \,}

\newcommand{\eqinlaw}{\, \stackrel{\mathcal{D}}{=} \,}
\newcommand{\bs}[1]{\boldsymbol{#1}}

\newcommand{\ssst}{\scriptscriptstyle}

\newcommand{\textfor}{\quad \text{for }}
\newcommand{\independent}{\perp \hspace{-2.3mm} \perp}

\newcommand{\BB}{\mathbb{B}}

\newcommand{\EE}{\mathbb{E}}

\newcommand{\NN}{\mathbb{N}}
\newcommand{\PP}{\mathbb{P}}
\newcommand{\QQ}{\mathbb{Q}}
\newcommand{\RR}{\mathbb{R}}

\newcommand{\Rplus}{\mathbb{R}_{+}}

\newcommand{\Zplus}{\mathbb{Z}_{+}}

\newcommand{\BBdot}{\dot{\BB}}

\newcommand{\linfty}{L_{\infty}}

\newcommand{\LebD}{\text{\textsl{Leb}}^D}
\newcommand{\dtv}{d_{TV}}

\newcommand{\dzero}{d_{0}}
\newcommand{\tdzero}{\td_{0}}
\newcommand{\done}{d_{1}}

\newcommand{\dtwo}{d_{2}}

\newcommand{\ftv}{\mathcal{F}_{TV}}

\newcommand{\ftwo}{\mathcal{F}_{2}}

\newcommand{\mvert}{\, \vert \,}
\newcommand{\abs}[1]{\lvert #1 \rvert}
\newcommand{\bigabs}[1]{\bigl| #1 \bigr|}
\newcommand{\Bigabs}[1]{\Bigl| #1 \Bigr|}
\newcommand{\biggabs}[1]{\biggl| #1 \biggr|}

\newcommand{\norm}[1]{\lVert #1 \rVert}
\newcommand{\bignorm}[1]{\bigl\lVert #1 \bigr\rVert}

\newcommand{\biggnorm}[1]{\biggl\lVert #1 \biggr\rVert}

\newcommand{\ldnorm}[1]{\norm{#1}_{L_{D}}}
\newcommand{\bigldnorm}[1]{\bignorm{#1}_{L_{D}}}

\newcommand{\linftynorm}[1]{\norm{#1}_{L_{\infty}}}

\newcommand{\bigglinftynorm}[1]{\biggnorm{#1}_{L_{\infty}}}

\newcommand{\indi}[1]{ 1_{#1} }

\newcommand{\varrhoint}{\varrho_{\textrm{int}}}
\newcommand{\varrhoext}{\varrho_{\textrm{ext}}}
\newcommand{\sigmaext}{\sigma_{\textrm{ext}}}
\newcommand{\etaext}{\eta_{\textrm{ext}}}

\newcommand{\Aint}{A_{\text{\textrm{int}}}}
\newcommand{\Aext}{A_{\text{\textrm{ext}}}}

\newcommand{\tb}{\tilde{b}}
\newcommand{\td}{\tilde{d}}
\newcommand{\tf}{\tilde{f}}
\newcommand{\tg}{\tilde{g}}
\newcommand{\tilh}{\tilde{h}}

\newcommand{\tp}{\tilde{p}}
\newcommand{\tr}{\tilde{r}}
\newcommand{\ts}{\tilde{s}}

\newcommand{\tx}{\tilde{x}}

\newcommand{\tC}{\widetilde{C}}

\newcommand{\tN}{\tilde{N}}

\newcommand{\tY}{\tilde{Y}}

\newcommand{\txi}{\tilde{\xi}}
\newcommand{\teta}{\tilde{\eta}}

\newcommand{\tmu}{\tilde{\mu}}

\newcommand{\tpi}{\tilde{\pi}}

\newcommand{\tvarrho}{\tilde{\varrho}}
\newcommand{\tsigma}{\tilde{\sigma}}

\newcommand{\mcb}{\mathcal{B}}

\newcommand{\bxc}{\BB_{\mcx}^{\text{\raisebox{0.5pt}{\hspace*{1pt}$c$}}}}
\newcommand{\bxpc}{\BB_{\text{\raisebox{0.5pt}[5pt][0pt]{$\mcx'$}}}^{\text{\raisebox{0.5pt}{\hspace*{1pt}$c$}}}}
\newcommand{\bxpcstd}{\BB_{\mcx'}^{\text{\raisebox{0.5pt}{\hspace*{1pt}$c$}}}}

\newcommand{\mcc}{\mathcal{C}}
\newcommand{\mce}{\mathcal{E}}

\newcommand{\mcg}{\mathcal{G}}

\newcommand{\mck}{\mathcal{K}}

\newcommand{\mcm}{\mathcal{M}}
\newcommand{\mcn}{\mathcal{N}}

\newcommand{\mcx}{\mathcal{X}}
\newcommand{\mcy}{\mathcal{Y}}
\newcommand{\msl}{\mathscr{L}}

\newcommand{\mfk}{\mathfrak{K}}
\newcommand{\mfm}{\mathfrak{M}}
\newcommand{\mfn}{\mathfrak{N}}
\newcommand{\mfp}{\mathfrak{P}}

\newcommand{\mcxpsq}{{\mcx'}^{\hspace*{0.02em}2}} 
\newcommand{\Sigmap}{\Sigma'}

\newcommand{\obenn}{^{(n)}}

\newcommand{\obenp}{^{{\ssst (} p {\ssst )}}}
\newcommand{\obenpi}{^{{\ssst (} \pi {\ssst )}}}

\newcommand{\tsum}{\textstyle{\sum}}

\newcommand{\RplusD}{\RR_{+}^D}

\newcommand{\RD}{\RR^D}
\newcommand{\cox}{\mathrm{Cox}}

\newcommand{\Deltaone}{\Delta_{\hspace*{-1pt}1}\hspace*{-0.9pt}}
\newcommand{\Deltatwo}{\Delta_{\hspace*{-0.5pt}2}\hspace*{-0.4pt}}

\newcommand{\ppth}[1]{#1_{\pi}}

\newcommand{\nin}{\noindent}

\newcommand{\kappatinv}{\kappa_T^{-1}}

\renewcommand{\theprop}{\arabic{section}.\Alph{prop}}

\begin{document}

\title{Distance estimates for dependent thinnings\\ of point processes with densities}
\author{By Dominic Schuhmacher\thanks{Work supported by the Swiss National Science Foundation under Grant
No.\ PBZH2-111668.}\\University of Western Australia} 
\maketitle

\begin{abstract}
   In [Schuhmacher, Electron.\ J.\ Probab.\ {\bf 10} (2005), 165--201] estimates of the Barbour-Brown
   distance $\dtwo$ between the distribution of a thinned point process and the distribution of a Poisson 
   process were derived by combining discretization with a result 
   based on Stein's method. In the present article we concentrate on point processes that have a density with
   respect to a Poisson process. For such processes we can apply a corresponding result directly without
   the detour of discretization and thus obtain better and more 
   natural bounds not only in $\dtwo$ but also in the stronger total variation metric. We give applications
   for thinning by covering with an independent Boolean model and ``Mat{\'e}rn type~I''-thinning of fairly
   general point processes. These applications give new insight into the respective models, and either
   generalize or improve earlier results.
\end{abstract} 


\section{Introduction}  \label{sec:intro}

We consider thinnings of simple point processes on a general compact metric space $\mcx$, where simple means that the probability of having multiple points at the same location is zero.
The thinning of such a process $\xi$ according to a $[0,1]$-valued measurable random field $\pi$ on $\mcx$ is the point process $\xi_{\pi}$, unique with regard to its distribution, that can be obtained in the following way: for any realizations $\xi(\omega)$ (a point measure on $\RD$) and $\pi(\omega,\cdot)$ (a function $\RD \to [0,1]$), look at each point $s$ of
$\xi(\omega)$ in turn, and retain it with probability $\pi(\omega,s)$, or delete it with
probability $1-\pi(\omega,s)$, independently of any retention/deletion decisions of other
points. Regard the points left over by this procedure as a realization of the thinned point process
$\xi_{\pi}$. For a formal definition see Section~\ref{sec:prelims}.
We usually refer to $\xi$ as \emph{the original process}, and to $\pi$ as \emph{the retention
field}. 

The following is a well-established fact: if we thin more and more, in the sense that we consider a sequence of retention fields $(\pi_n)_{n \in \NN}$ with \raisebox{0pt}[11.5pt][0pt]{$\sup_{x \in \mcx} \pi_n(x) \inlawto 0$} as $n \to \infty$, and compensate for the thinning by choosing point processes $\xi_n$ whose intensity increases as $n$ goes to infinity in a way that is compatible with the retention fields, then we obtain convergence towards a Cox process. The theorem below was shown in \cite{kallenberg75} for constant deterministic $\pi_n$, and generalized in \cite{brown79} to general $\pi_n$. In order to specify what choice of the sequence $(\xi_n)$ is compatible with $(\pi_n)$, we introduce the random measure $\Lambda_n$ that is given by $\Lambda_n(A) := \int_{A} \pi_n(x) \; \xi_n(dx)$ for every Borel set $A$ in $\mcx$. Convergence of random measures, and in particular of point processes, is defined via the convergence of expectations of bounded continuous functions, where continuity is in terms of the vague topology (for details see \cite{kallenberg86}, Section 4.1).
\begin{thm}[Kallenberg, Brown] \label{thm:kb}
   For the sequences $(\pi_n)_{n \in \NN}$ and $(\xi_n)_{n \in \NN}$ introduced above, we obtain
   convergence in distribution of the thinned sequence $\bigl( (\xi_n)_{\pi_n} \bigr)_{n \in \NN}$ towards a
   point process $\eta$ if and only if there is a random measure $\Lambda$ on $\mcx$ such that
   \raisebox{0pt}[10.75pt][0pt]{$\Lambda_n \inlawto \Lambda$} as $n \to \infty$.
   In this case $\eta \sim \cox(\Lambda)$, i.e.\ $\eta$ is a Cox process with directing measure $\Lambda$. 
\end{thm}

In \cite{schumi2} the above setting was considered for the situation that the $\xi_n$ are point processes on $[0,1]^D$ that are obtained from a single point process $\xi$ by gradual contraction of $\RplusD$ using the functions $\kappa_n$ given by $\kappa_n(x) := (1/n)x$ for every $x \in \RplusD$ (for notational convenience in the proofs, the order of contracting and thinning was interchanged). Under the additional assumption that $\xi$ and $\pi_n$ satisfy mixing conditions, which makes it plausible for the limiting process in Theorem~\ref{thm:kb} to be Poisson (see the remark after Theorem~1.A of \cite{schumi2}), several upper bounds for the Wasserstein-type distance $\dtwo$ between the distribution of the thinned and
contracted process and a suitable Poisson process distribution
were obtained under various conditions. These results were derived by discretizing the thinned process and the limiting Poisson process, and applying then a discrete version (essentially formulated as Theorem~10.F in \cite{bhj92}) of the ``local Barbour-Brown theorem'', Theorem~3.6 in \cite{bb92}, which is based on Stein's method (see \cite{stein72} and~\cite{steinintro}). 

Although the bounds were of good quality and have proved their usefulness in several applications, they had some shortcomings, which were mainly related to the fact that they were still expressed in terms of discretization cuboids. This made the results rather unpleasant to apply in many situations where truly non-discrete point processes were considered.

The present article deals with these issues by offering a proof that makes direct use of Theorem~3.6 in \cite{bb92} without the intermediate step of discretization. There are several advantages of this approach: the proofs become simpler and more elegant; the upper bounds are much more natural, more easily applied, hold for the most part under more general conditions and are qualitatively slightly better; what is more, we can obtain a bound not only for the $\dtwo$-distance, but also for the stronger total variation distance. In order to apply Theorem~3.6, we restrict ourselves to point processes $\xi$ that have a density with respect to the distribution of a simple Poisson process, which is a natural and common choice for a reference measure and leaves us with a very rich class of processes.
The simplicity of the Poisson process implies the simplicity of $\xi$. It should be noted, however, that any non-simple point process can be turned into a simple one by lifting it to a bigger space $\mcx \times K$. This space has to be compact in our setting, so that e.g.\ $K=\{1,2, \ldots, k\}$ (if there is a maximal number $k$ of points per location), $K=[0,1]$, or $K=[0,\eps]$ for $\eps$ small are typical choices.

We start out in Section~\ref{sec:prelims} by giving an overview of the technical background and some notation needed to formulate and prove the main results. These results are then presented in Section~3. We provide
upper bounds for the $\dtv$- and the $\dtwo$-distances between $\msl( \xi_\pi )$ and a suitable Poisson process distribution, first in a general setting, and then for a number of important special cases. The last of these special cases (see Corollary~\ref{cor:maincontracted}) is suitable for comparison with the upper bounds in \cite{schumi2}. Finally, in Section~\ref{sec:applications}, two ``extreme'' applications of the main results are studied, in which the retention field takes only the values $0$ and $q \in [0,1]$ and is either independent or completely determined by $\xi$. In the first application a point is always deleted if covered by a certain independent Boolean model, and is retained with probability $q$ otherwise. In the second application, a point is always deleted if there is any other point present within a fixed distance $r$ (following the construction of the Mat{\'e}rn type I hard core process), and again is retained with probability $q$ otherwise.

\section{Preliminaries} \label{sec:prelims}

We first introduce some basic notation and conventions, before giving an overview of some of the theoretical background and presenting the more technical definitions in the various subsections.

Let $(\mcx, \dzero)$ always be a compact metric space with $\dzero \leq 1$ that admits a finite diffuse measure
$\alpha \neq 0$, where diffuse means that $\alpha(\{x\}) = 0$ for every $x \in \mcx$.
Denote by $\mcb$ the Borel $\sigma$-algebra on $\mcx$, and by $\mcb_A$ the trace $\sigma$-algebra $\mcb \vert_A = \{ B \cap A; B \in \mcb \}$ for any set $A \subset \mcx$. Furthermore, write $\mfm$ for the
space of finite measures on $\mcx$, and equip it with the vague topology (see \cite{kallenberg86},
Section~15.7) and the corresponding Borel $\sigma$-algebra $\mcm$ (see \cite{kallenberg86}, Lemma~4.1 and Section~1.1). Do the same for the
subspace $\mfn \subset \mfm$ of finite point measures, and denote the corresponding
$\sigma$-algebra by~$\mcn$. Write furthermore $\mfn^{*} := \{ \varrho \in \mfn; \, \varrho(\{x\}) \leq 1 \text{ for all } x \in \mcx \}$ for the $\mcn$-measurable set of simple point measures. A \emph{random measure} on $\mcx$ is a random element of~$\mfm$, and a \emph{point process} on $\mcx$ a random element of $\mfn$. A point process $\xi$ is called \emph{simple} if $\PP[\xi \in \mfn^{*}] = 1$. By $\poisson(\lambda)$ we denote the distribution of the Poisson process on $\mcx$ with
intensity measure $\lambda$ if $\lambda \in \mfm$, and the Poisson distribution with parameter
$\lambda$ if $\lambda$ is a positive real number. 

We think of measures (random or otherwise) always as being defined on all of $\mcx$. Thus, for any measure $\mu$ on $\mcx$ and any $A \in \mcb$, we denote by $\mu \vert_A$ the measure on $\mcx$ that is given by $\mu \vert_A(B) := \mu(A \cap B)$ for all $B \in \mcb$.
Let $\mfm(A) := \{ \mu \vert_A ; \, \mu \in \mfm \}$, $\mfn(A) := \{ \varrho \vert_A ; \, \varrho \in \mfn \}$, $\mcm(A) := \mcm \vert_{\mfm(A)} = \{ C \cap \mfm(A) ; \, C \in \mcm \}$, and $\mcn(A) := \mcn \vert_{\mfn(A)} = \{ C \cap \mfn(A) ; \, C \in \mcn \}$. Furthermore, set $\mfn^{*}(A) := \mfn(A) \cap \mfn^{*}$. Sometimes absolute value bars are used to denote the total mass of a measure, i.e.\ $\abs{\mu} := \mu(\mcx)$ for
any $\mu \in \mfm$.

For $\sigma \in \mfn^{*}$, we do not notationally distinguish between the point measure and its support. Like that we can avoid having to enumerate the points of the measure, which sometimes saves us from tedious notation. Thus, the notation $\int f \: d\sigma = \sum_{s \in \sigma} f(s)$ may be used instead of writing $\int f \: d\sigma = \sum_{i = 1}^v f(s_i)$, where $\sigma = \sum_{i=1}^v \delta_{s_i}$. The same concept is extended to the realizations of a simple point process, as long as it is not important what happens on the null set of point patterns with multiple points. Hence we may also write $\EE \bigl( \int f \: d\xi \bigr) = \EE \bigl( \sum_{s \in \xi} f(s) \bigr)$ if $\xi$ is simple. In order to facilitate the reading of certain formulae, we furthermore make the convention of using the letters $x, \tx, y$ for general elements of the state space $\mcx$, while $s, \ts, t$ are reserved for points of a point pattern in~$\mcx$. Finally, we sometimes omit the addition ``almost surely'' or ``(for) almost every \ldots'' for equations and inequalities between functions on measure spaces if it is evident and of no importance that the relation does not hold pointwise.

\subsection{Densities with respect to the standard Poisson process distribution $\bs{P_1}$}

Let $\alpha \neq 0$ be a fixed diffuse measure in $\mfm$,
which we will regard as our reference measure
on~$\mcx$. Typically, if (a superset of) $\mcx$ has a suitable group structure, $\alpha$ is chosen to be (the restriction of) the Haar measure.
If $\mcx \subset \RD$, we usually choose $\alpha = \LebD \vert_{\mcx}$ of course. We write $P_1 := \poisson(\alpha)$ for the distribution of what we call the \emph{standard Poisson process} on $\mcx$, and $P_{1,A} := \poisson(\alpha \vert_A)$ for the distribution of the Poisson process on $\mcx$ with expectation measure $\alpha \vert_A$. It is convenient to admit also $\alpha(A) = 0$, in which case $P_{1,A} = \delta_0$, where $0$ denotes the zero measure on $\mcx$.

A popular way of specifying a point process distribution is by giving its Radon-Nikodym density with respect to the distribution of the standard Poisson process (if the density exists; see \cite{moellerwaage04}, Sections~6.1 and~6.2 for a number of examples). The following simple lemma will be useful.
\begin{lem} \label{lem:hompoisdensity}
   For any $a > 0$, a density $f_A$ of $P_{a,A} := \poisson(a \! \cdot \! \alpha \vert_{A})$ with respect to
   $P_{1,A}$ is given by
   \begin{equation*}
      f_A(\varrho) = e^{(1-a) \alpha(A)} a^{\abs{\varrho}}
   \end{equation*}
   for every $\varrho \in \mfn(A)$.
\end{lem}
\begin{proof}
   See \cite{moellerwaage04}, Proposition~3.8(ii), for the case $\mcx \subset \RD$, $\alpha = \LebD \vert_
   {\mcx}$, and a more general Poisson process. The same proof holds for general $\mcx$ and $\alpha$.
\end{proof}

To avoid certain technical problems, we require our densities to be hereditary whenever we are dealing with conditional densities. 
\begin{defi}
   A function $f: \mfn \to \Rplus$ is called \emph{hereditary} if $f(\varrho) = 0$ implies $f(\sigma) = 0$
   whenever $\varrho, \sigma \in \mfn$ with $\varrho \subset \sigma$.
\end{defi} 
\nin The point processes on $\mcx$ that have a hereditary density with respect to $P_1$ form the class of \emph{Gibbs point processes}. These include pairwise and higher order interaction processes (see \cite{moellerwaage04}, Section~6.2). The general form of a Gibbs process density is given in Definition~4.2 of \cite{baddeley07}.

\subsection{Thinnings} \label{ssec:thinnings}

In what follows, let $\xi$ always be a point process on $\mcx$ that has density $f$ with respect to $P_1$, and let $\pi := ( \pi(\cdot, x); x \in \mcx )$ be a
$[0,1]$-valued random field on $\mcx$ that is measurable in the sense that the mapping $\Omega
\times \mcx \to [0,1]$, $(\omega,x) \mapsto \pi(\omega,x)$ is $\sigma(\pi) \otimes
\mcb$-measurable.
In the main part we strengthen the technical conditions on $\xi$ and $\pi$ somewhat in order to avoid
some tedious detours in the proofs.

We use the definition from \cite{schumi2} for the $\pi$-thinning of $\xi$.
\begin{defi}[Thinning]
First, assume that $\xi =
\sigma = \sum_{i=1}^v \delta_{s_i}$ and $\pi=p$ are non-random, where $v \in \Zplus$, $s_i \in \mcx$, and $p$ is a measurable function $\mcx \to [0,1]$. Then, a $\pi$-thinning of
$\xi$ is defined as $\xi_{\pi} = \sum_{i=1}^v X_i \delta_{s_i}$, where the $X_i$ are
independent indicator random
variables with expectations $p(s_i)$, respectively. Under these circumstances, 
$\xi_{\pi}$ has a distribution $P_{(\sigma, \, p)}$ that does not depend on the chosen
enumeration of~$\sigma$. We obtain the general $\pi$-thinning from this by randomization,
that is by the condition $\PP [ \ppth{\xi} \in \cdot \mvert \xi, \pi ] = P_{(\xi, \pi)}$ (it is straightforward to see that $P_{(\xi, \pi)}$ is a $\sigma(\xi,
\pi)$-measurable family in the sense that $P_{(\xi, \pi)}(D)$ is $\sigma(\xi, \pi)$-measurable
for every $D \in \mfn$). Note that the distribution of $\ppth{\xi}$ is uniquely determined by
this procedure.
\end{defi}

The following lemma gives the first two factorial moment measures of $\xi_{\pi}$.
For $\varrho = \sum_{i=1}^v \delta_{s_i} \in \mfn$, write $\varrho^{[2]} := \sum_{i,j=1, i \neq j}^v
\delta_{(s_i,s_j)}$ for the factorial product measure on $\mcx \times \mcx$. Remember that the expectation measure $\mu_1$ of $\xi$ is defined by $\mu_1(A) := \EE(\xi(A))$ for every $A \in \mcb$, and that the second factorial moment measure $\mu_{[2]}$ of $\xi$ is defined to be the expectation measure of $\xi^{[2]}$. Let $\Lambda$ be the random measure on $\mcx$ that is given by $\Lambda(A) := \int_A \pi(x) \; \xi(dx)$ for $A \in \mcb$ (cf.\ $\Lambda_n$ in Section~\ref{sec:intro}).
\begin{lem} \label{lem:mms}
   We obtain for the expectation measure $\mu_1\obenpi$ and the second factorial moment measure
   $\mu_{[2]}\obenpi$ of $\xi_{\pi}$
   \begin{enumerate}
      \item ${\displaystyle \mu_1\obenpi(A) = \EE \biggl( \int_A \pi(x) \; \xi(dx) \biggr) = \EE
         \bigl( \Lambda(A) \bigr)}$ for every $A \in \mcb$;
      \item ${\displaystyle \mu_{[2]}\obenpi(B) = \EE \biggl( \int_B \pi(x) \pi(\tx) \;   
         \xi^{[2]} \bigl( d(x,\tx) \bigr) \biggr) =: \EE \bigl( \Lambda^{[2]}(B) \bigr)}$ for every
         $B \in \mcb^2$.
   \end{enumerate}
\end{lem}
\begin{proof}
Write $\xi = \sum_{i=1}^V \delta_{S_i}$, where $V$ and $S_i$ are $\sigma(\xi)$-measurable random elements
with values in $\Zplus$ and $\mcx$, respectively. Such a representation exists by Lemma~2.3 in \cite{kallenberg86}.
\begin{enumerate}
   \item For every $A \in \mcb$ we have
   \begin{equation*}
   \begin{split}
      \mu_1\obenpi(A) = \EE \bigl( \xi_{\pi}(A) \bigr) &= \EE \bigl( \EE \bigl( \xi_{\pi}(A) \bigm| \xi,
         \pi \bigr) \bigr) \\
      &= \EE \biggl( \sum_{i=1}^V \pi(S_i) \hspace*{1pt} {\rm I}[S_i \in A] \biggr)
         = \EE \biggl( \int_A \pi(x) \; \xi(dx) \biggr).
   \end{split}
   \end{equation*}
   \item For every $B \in \mcb^2$ we have
   \begin{equation*}
   \begin{split}
      \mu_{[2]}\obenpi(B) = \EE \bigl( \xi_{\pi}^{[2]}(B) \bigr) &= \EE \bigl( \EE \bigl(
         \xi_{\pi}^{[2]}(B) \bigm| \xi, \pi \bigr) \bigr) \\
      &= \EE \Biggl( \sum_{\substack{i, j=1 \\ i \neq j}}^V \pi(S_i) \pi(S_j) \hspace*{1pt} {\rm I}[(S_i, 
         S_j) \in B] \Biggr)
      = \EE \biggl( \int_B \pi(x) \pi(\tx) \; \xi^{[2]}\bigl(d(x,\tx)\bigr) \biggr).
   \end{split}
   \end{equation*}
\end{enumerate}
\end{proof}

\subsection{Metrics used on the space of point process distributions}

We use two metrics on the space $\mfp(\mfn)$ of probability measures on $\mfn$. The one that is more widely known is the total variation metric, which can be defined on any space of probability measures. For $P, Q \in \mfp(\mfn)$ it is given by
\begin{equation*}
   \dtv(P,Q) = \sup_{C \in \mcn} \bigabs{P(C) - Q(C)}
\end{equation*}
or, equivalently, by
\begin{equation} \label{eq:dtvcoupling}
   \dtv(P,Q) = \min_{\substack{\xi_1 \sim P \\ \xi_2 \sim Q}} \, \PP[ \xi_1 \neq \xi_2]. 
\end{equation}
See \cite{bhj92}, Appendix A.1, for this and other general results about the total variation metric.

The second metric we use is a Wasserstein type of metric introduced by Barbour and Brown in \cite{bb92}, and denoted by $\dtwo$. It is often a more natural metric to use on $\mfp(\mfn)$ than $\dtv$, because it takes the metric $\dzero$ on $\mcx$ into account and metrizes convergence in distribution of point processes. The total variation metric on the other hand is strictly stronger, and at times too strong to be useful.

Define the
\emph{$\done$-distance} between two point measures $\varrho_1 = \sum_{i=1}^{\abs{\varrho_1}} \delta_{s_{1,i}}$ and $\varrho_2 = \sum_{i=1}^{\abs{\varrho_2}} \delta_{s_{2,i}}$ in $\mfn$ as
\begin{equation}
   \done(\varrho_1, \varrho_2) := \begin{cases} 1 &\text{if }
   \abs{\varrho_1} \neq \abs{\varrho_2},  \\ 
   \min_{\tau \in \Sigma_v} \bigl[ \frac{1}{v} \tsum_{i=1}^v 
         \dzero(s_{1,i}, s_{2,\tau(i)}) \bigr] &\text{if }
   \abs{\varrho_1} = \abs{\varrho_2} = v > 0,  \\
   0 &\text{if } \abs{\varrho_1} = \abs{\varrho_2} = 0, \end{cases}
\end{equation}
where $\Sigma_v$ denotes the set of permutations of $\{ 1, 2, \ldots, v \}$. 
It can be seen that $(\mfn,\done)$ is a complete, separable metric space and that $\done$ is bounded by 1.
Now let $\ftwo
:= \{ f : \mfn \to \RR; \; \lvert f(\varrho_1) - f(\varrho_2) \rvert
\leq \done(\varrho_1,\varrho_2) \text{ for all $\varrho_1, \varrho_2 \in \mfn$} \}$, and define the
\emph{$\dtwo$-distance} between two measures $P, Q \in \mfp(\mfn)$ as
\begin{equation*}
      \dtwo(P, Q) := \sup_{f \in \ftwo} \left| \int f \; dP - \int 
      f \; dQ \right|
\end{equation*}
or, equivalently, as
\begin{equation}  \label{eq:krdtwo}
      \dtwo(P,Q) = \min_{\substack{\xi_1 \sim P \\ \xi_2 \sim Q}} \, 
      \mathbb{E} \done(\xi_1, \xi_2).
\end{equation}
See \cite{schumi2}, \cite{mythesis}, and \cite{xia05} for this and many other results about the Barbour-Brown metric $\dtwo$. By $\eqref{eq:dtvcoupling}$ and $\eqref{eq:krdtwo}$ we obtain that $\dtwo \leq \dtv$.

\subsection{Distance estimates for Poisson process approximation of point processes with a spatial dependence structure} \label{ssec:bb}

In this subsection, a theorem is presented that provides upper bounds for the total variation and
$\dtwo$ distances between the distribution of a point process with a spatial dependence structure and a suitable Poisson process distribution. This result is very similar to Theorems~2.4 and~3.6 in \cite{bb92}, but deviates in several minor aspects, two of which are more pronounced: first, we use densities with respect to a Poisson process distribution instead of Janossy densities, which simplifies certain notations considerably; secondly, we take a somewhat different approach for controlling the long range dependence within $\xi$ (see the term for $\bphi(x)$ in Equation~\eqref{eq:bphinonher}), which avoids the imposition of an unwelcome technical condition (cf.\ Remark~\ref{rem:bbamend}).

Let $\xi$ be a point process on $\mcx$ whose distribution has a density $f$ with respect to $P_1$ and whose expectation measure $\mu = \mu_1$ is finite. Then $\mu$ has a density $u$ with respect to $\alpha$ that is given by
\begin{equation} \label{eq:u}
   u(x) = \int_{\mfn} f(\varrho + \delta_x) \; P_1(d \varrho)
\end{equation}
for $\alpha$-almost every $x \in \mcx$, which is obtained as a special case of
Equation~\eqref{eq:pregnzconsequence} below by choosing $N_x := \{x\}$.

For $A \in \mcb$, we set
\begin{equation} \label{eq:densityofrestricted}
   f_A(\varrho) := \int_{\mfn(A^c)} f(\varrho + \tvarrho) \; P_{1,A^c}(d \tvarrho)
\end{equation}
for every $\varrho \in \mfn(A)$,
which gives a density $f_A: \mfn \to \Rplus$ of the distribution of $\xi \vert_A$ with respect to $P_{1,A}$ (we extend $f_A$ on $\mfn \setminus \mfn(A)$ by setting it to zero). 
This can be seen by the fact that integrating $f_A$ over an arbitrary set $C \in \mfn(A)$ yields
\begin{equation*}
\begin{split}
   \int_C f_A(\varrho) \; P_{1,A}(d\varrho) &= \int_{\mfn(A)} \int_{\mfn(A^c)} {\rm I}[ \varrho \in 
      C] f(\varrho + \tvarrho) \; P_{1,A^c}(d \tvarrho) \; P_{1,A}(d \varrho) \\
   &= \int_{\mfn} {\rm I} [\sigma \vert_A \in C] f(\sigma) \; P_1(d \sigma) \\
   &= \PP[\xi \vert_A \in C],
\end{split}
\end{equation*}
where we used that $( \eta \vert_A, \eta \vert_{A^c}) \sim P_{1,A} \otimes P_{1,A^c}$ for $\eta \sim P_1$. Note that the argument remains correct if either $\alpha(A)$ or $\alpha(A^c)$ is zero.

We introduce a \emph{neighborhood structure} $(N_x)_{x \in \mcx}$ on $\mcx$, by which we mean any collection of sets that satisfy $x \in N_x$ for every $x \in \mcx$ (note that we do not assume $N_x$ to be a neighborhood of $x$ in the topological sense). We require this neighborhood structure to be measurable in the sense that
\begin{equation} \label{eq:nhcondition}
   N(\mcx) := \{ (x,y) \in \mcx^2; y \in N_x \} \in \mcb^2.
\end{equation}
This measurability condition is slightly stronger than the ones required in \cite{bb92} (see the discussion before Remark~2.1  in \cite{chenxia04} for details), but quite a bit more convenient. The $N_x$ play the role of regions of strong dependence: it is advantageous in view of Theorem~\ref{thm:bb} below to choose $N_x$ not too large, but in such a way that the point process $\xi$ around the location $x$ depends only weakly on $\xi \vert_{N_x^c}$. Write $\dot{N}_x$ for $N_x \setminus \{x\}$.

We use the following crucial formula about point process densities, which is proved as Proposition~\ref{prop:pregnz} in the appendix. For any
non-negative or bounded measurable function $h : \mcx \times \mfn \to~\RR$, we have
\begin{equation} \label{eq:pregnz}
   \EE \biggl( \int_{\mcx} h( x, \xi \vert_{N_x^c}) \; \xi(dx) \biggr) = \int_{\mcx} \int_{\mfn(N_x^c)} h
   (x,\varrho) f_{N_x^c \cup \{x\}}(\varrho + \delta_x) \; P_{1,N_x^c}(d \varrho) \; \alpha(dx).
\end{equation} 
As an important consequence we obtain by choosing $h(x,\varrho) := 1_A(x)$ that
\begin{equation} \label{eq:pregnzconsequence}
   u(x) = \int_{\mfn(N_x^c)} f_{N_x^c \cup \{x\}} (\varrho + \delta_x) \; P_{1, N_x^c}(d \varrho)
\end{equation}
for $\alpha$-almost every $x \in \mcx$.

In many of the more concrete models, the density $f$ of $\xi$ is hereditary and therefore $\xi$ is a Gibbs process, in which case the above expressions can be simplified by introducing conditional densities. 
Let $\mfk := \bigcup_{x \in \mcx} \bigl( \{x\} \times \mfn(N_x^c) \bigr) \subset \mcx \times \mfn$ and define a mapping $g: \mcx \times \mfn \to \Rplus$ by
\begin{equation} \label{eq:gdefi}
   g(x; \varrho) := \frac{f_{N_x^c \cup \{x\}}(\varrho + \delta_x)}{f_{N_x^c}(\varrho)}
\end{equation}
for $(x, \varrho) \in \mfk$ and $g(x; \varrho) := 0$ otherwise, where the fraction in \eqref{eq:gdefi} is taken to be zero if the denominator (and hence by hereditarity also the enumerator) is zero.
For $(x,\varrho) \in \mfk$ the term $g(x;\varrho)$ can be interpreted as the conditional density of a point at $x$ given the configuration of $\xi$ outside of $N_x$ is $\varrho$.
Equation~\eqref{eq:pregnz} can then be replaced by the following result, which is a generalization of the Nguyen-Zessin Formula (see \cite{nguyenzessin79}, Equation~(3.3)): for any non-negative or bounded measurable function $h : \mcx \times \mfn \to \RR$, we have
\begin{equation} \label{eq:nzplus}
   \EE \biggl( \int_{\mcx} h( x, \xi \vert_{N_x^c}) \; \xi(dx) \biggr) = \int_{\mcx} \EE \bigl( h( x, \xi
   \vert_{N_x^c}) g(x; \xi \vert_{N_x^c}) \bigr) \; \alpha(dx). 
\end{equation}
This formula was already stated as Equation~(2.7) in~\cite{bb92} for functions $h$ that are constant in $x$ 
and as Equation~(2.10) in \cite{chenxia04} for general functions, both times however under too wide conditions. See Corollary~\ref{cor:gnz} for the proof and Remark~\ref{rem:bbamend} for an example that shows that an additional assumption, such as hereditarity, is required. As an analog to $\eqref{eq:pregnzconsequence}$, we obtain
\begin{equation*}
   u(x) = \EE \bigl( g(x; \xi \vert_{N_x^c}) \bigr)
\end{equation*}
for $\alpha$-almost every $x \in \mcx$.

We are now in a position to derive the required distance bounds.
\begin{thm}[based on Barbour and Brown \cite{bb92}, Theorems 2.4 and 3.6] \label{thm:bb}
   Suppose that $\xi$ is a point process which has density $f$ with respect to $P_1$ and finite expectation 
   measure~$\mu$. Let furthermore $(N_x)_{x \in \mcx}$ be a neighborhood structure that is measurable in 
   the sense of Condition~\eqref{eq:nhcondition}. Then
   \begin{equation*}
   \begin{split}
      \text{(i) } \, &\dtv \bigl( \msl(\xi), \poisson(\mu) \bigr) \, \leq \, \int_{\mcx} \mu(N_x) \;
         \mu(dx) + 
         \EE \biggl( \int_{\mcx} \xi(\dot{N}_x) \; \xi(dx) \biggr) + \int_{\mcx} \bphi(x) \; \alpha(dx) ;   
         \\[2mm]
      \text{(ii) } \, &\dtwo \bigl( \msl(\xi), \poisson(\mu) \bigr) \, \leq \, M_2(\mu) \biggl[ \int_{\mcx} 
         \mu(N_x) \; \mu(dx) + \EE \biggl( \int_{\mcx} \xi(\dot{N}_x) \; \xi(dx) \biggr) \biggr] 
          + M_1(\mu) \hspace*{-0.5mm} \int_{\mcx} \bphi(x) \; \alpha(dx) ;
   \end{split}
   \end{equation*}
   where
   \begin{equation*}
      M_1(\mu) = \min \biggl( 1, \frac{1.647}{\sqrt{\abs{\mu}}} \biggr), \:
      M_2(\mu) = \min \biggr( 1, \frac{11}{6 \abs{\mu}} \Bigl( 1 + 2 \log^{+} \Bigl(
         \frac{6 \abs{\mu}}{11} \Bigr) \Bigr) \biggr),
   \end{equation*}
   and 
   \begin{equation} \label{eq:bphinonher}
   \begin{split}
      \bphi(x) &= \int_{\mfn(N_x^c)} \bigabs{ f_{N_x^c \cup \{x\}}(\varrho +
         \delta_x) - f_{N_x^c}(\varrho) u(x) } \; P_{1,N_x^c}(d \varrho) \\
      &= 2 \sup_{C \in \mcn(N_x^c)}
         \biggabs{\int_C \bigl[ f_{N_x^c \cup \{x\}}(\varrho +
         \delta_x) - f_{N_x^c}(\varrho) u(x) \bigr] \; P_{1,N_x^c}(d \varrho)}.
   \end{split}
   \end{equation}
   If $f$ is hereditary, $\bphi$ can be rewritten as
   \begin{equation} \label{eq:bphiher}
      \bphi(x) = \EE \bigabs{ g(x; \xi \vert_{N_x^c}) - u(x) }. 
   \end{equation}
\end{thm}
\begin{rem} \label{rem:wdt}
   We refer to the three summands in the upper bound of Theorem 2.C.(i) as basic term, strong dependence 
   term, and weak dependence term, respectively. The basic term depends only on the
   first order properties of $\xi$ and on $\alpha(N_x)$. The strong dependence term controls what    
   happens within the neighborhoods of strong dependence and is small if $\alpha(N_x)$ is not too big 
   and if there is not too much positive short range dependence within $\xi$. Finally, the weak dependence term 
   is small if there is only little long range dependence.
\end{rem}
\begin{rem}
   Theorem~5.27 in \cite{xia05} (which is based on several earlier results by various authors) gives an 
   alternative upper bound for the $\dtwo$-distance above, which when carefully further estimated is in many
   situations superior to the one in \cite{bb92}, insofar as the logarithmic term in 
   $M_2(\mu)$ can often be disposed of.
   After applying the same modifications as in the proof of Theorem~\ref{thm:bb} below,
   it can be seen that
   \begin{equation*}
   \begin{split}
      \dtwo \bigl( \msl(\xi), \poisson(\mu) \bigr) \, &\leq \, \int_{\mcx} 
         \EE \Bigl( \Bigl( \frac{3.5}{\abs{\mu}} + \frac{2.5}{\xi(N_x^c)+1} \Bigr) \xi(N_x) \Bigr) \;
         \mu(dx) \\
      &\hspace*{4mm} + \EE \biggl( \int_{\mcx} \Bigl( \frac{3.5}{\abs{\mu}} + \frac{2.5}{\xi(N_x^c)+1} \Bigr)
         \xi(\dot{N}_x) \; \xi(dx) \biggr)
         + M_1(\mu) \hspace*{-0.5mm} \int_{\mcx} \bphi(x) \; \alpha(dx).
   \end{split}
   \end{equation*} 
   Since working with this inequality requires a more specialized treatment of the thinnings in our main 
   result, and since the benefit of removing the 
   logarithmic term is negligible for most practical purposes, we do not use this bound in the present    
   article.
\end{rem}
\begin{proof}[Proof of Theorem~\ref{thm:bb}]
   Following the proof of Theorems~2.4 and~3.6 in \cite{bb92} (applying Equations~(2.9) and~(2.10) of
   \cite{bb92}, but not Equation~(2.11)), we obtain by using Stein's method that
   \begin{equation} \label{eq:bborig}
   \begin{split}
      \bigabs{\EE \tf(\xi) - \EE \tf(\eta)} &\leq \Deltatwo \tilh \, \biggl[ \int_{\mcx} \EE \xi
         (N_x) \; \mu(dx) + \EE \biggl( \int_{\mcx} \xi(\dot{N}_x) \; \xi(dx) \biggr) \biggr] \\
      &\hspace*{9mm} {} + \biggabs{\EE \int_{\mcx} \bigl[ \tilh(\xi \vert_{N_x^c} + \delta_x) - 
         \tilh(\xi \vert_{N_x^c}) \bigr] \bigl( \xi(dx) - \mu(dx) \bigr)},
   \end{split}
   \end{equation}
   for every $\tf \in \ftv = \{ \indi{C}; \, C \in \mcn \}$ [or $\tf \in \ftwo$ in the case of
   statement~(ii)], where $\eta \sim \poisson(\mu)$ and $\tilh := \tilh_{\tf}$ are the solutions of the
   so-called Stein equation (see \cite{bb92}, Equation~(2.2)), which have maximal first and second differences 
   \begin{equation*}
      \Deltaone \tilh := \sup_{\varrho \in \mfn, \, x \in \mcx} \bigabs{\tilh(\varrho + \delta_x) -
         \tilh(\varrho)} 
   \end{equation*}
   and
   \begin{equation*}
      \Deltatwo \tilh := \sup_{\varrho \in \mfn, \, x,y \in \mcx} \bigabs{\tilh(\varrho + \delta_x + 
      \delta_y) - \tilh(\varrho + \delta_x) - \tilh(\varrho + \delta_y) + \tilh(\varrho)}
   \end{equation*}
   that are bounded by $1$ [or $\Deltaone \tilh \leq M_1(\mu)$ and $\Deltatwo \tilh \leq M_2(\mu)$ in the
   case of statement~(ii); see \cite{xia05}, Propositions 5.16 and 5.17].
   
   All that is left to do is bounding the term in the second line of \eqref{eq:bborig}, which is done as 
   follows. 
   Setting $C_x^{+} := \bigl\{ \varrho \in \mfn(N_x^c) ; \, f_{N_x^c \cup \{x\}}(\varrho + \delta_x) >
   f_{N_x^c}(\varrho) u(x) \bigr\} \in \mcn(N_x^c)$, we obtain
   \begin{equation} \label{eq:bbnew}
   \begin{split}
      \biggabs{ \EE \int_{\mcx} \bigl[ \tilh(\xi &\vert_{N_x^c} + \delta_x) - \tilh(\xi
         \vert_{N_x^c}) \bigr] \bigl( \xi(dx) - \mu(dx) \bigr) }\\
      &= \biggabs{ \int_{\mcx} \int_{\mfn(N_x^c)} \bigl[ \tilh(\varrho + \delta_x) - \tilh(\varrho)
         \bigr] \bigl(
         f_{N_x^c \cup \{x\}}(\varrho + \delta_x) - f_{N_x^c}(\varrho) u(x) \bigr) \; P_{1,N_x^c}(d  
         \varrho) \; \alpha(dx) } \\
      &\leq \Deltaone \tilh \int_{\mcx} \int_{\mfn(N_x^c)} \bigabs{ f_{N_x^c \cup \{x\}}(\varrho + 
         \delta_x) - f_{N_x^c}(\varrho) u(x) } \; P_{1,N_x^c}(d \varrho) \; \alpha(dx) \\
      &= 2 \, \Deltaone \tilh \int_{\mcx} \int_{C_x^{+}} \bigl( f_{N_x^c \cup \{x\}}(\varrho + \delta_x) -
         f_{N_x^c}(\varrho) u(x) \bigr) \; P_{1,N_x^c}(d \varrho) \; \alpha(dx) \\
      &= 2 \, \Deltaone \tilh \int_{\mcx} \sup_{C \in \mcn(N_x^c)} \biggabs{
         \int_{C} \bigl( f_{N_x^c \cup \{x\}}(\varrho + \delta_x) -
         f_{N_x^c}(\varrho) u(x) \bigr) \; P_{1,N_x^c}(d \varrho) } \; \alpha(dx),
   \end{split}
   \end{equation}
   where we use Equation~\eqref{eq:pregnz}
   for the second line and 
   \begin{equation*}
      \int_{\mfn(N_x^c)} \bigl( f_{N_x^c \cup \{x\}}(\varrho + \delta_x) - f_{N_x^c}(\varrho) u(x) \bigr)  
      \; P_{1,N_x^c}(d \varrho) = 0
   \end{equation*}
   for the fourth line, which follows from Equation~\eqref{eq:pregnzconsequence}. The integrands with
   respect to $\alpha(dx)$ in the last three lines of \eqref{eq:bbnew} are all equal to
   \begin{equation*}
      \int_{\mfn} \bigabs{ f_{N_x^c \cup \{x\}}(\sigma \vert_{N_x^c} + \delta_x) -
      f_{N_x^c}(\sigma \vert_{N_x^c} ) u(x) } \; P_{1}(d \sigma) 
   \end{equation*}
   and hence $\mcb$-measurable by the defintion of $f_A$ and the fact that Condition~\eqref{eq:nhcondition}   
   implies the measurability of the mapping $\bigl[ \mcx \times \mfn \to \mcx \times \mfn, (x,\sigma) \mapsto 
   (x, \sigma \vert_{N_x^c}) \bigr]$ (see \cite{chenxia04}, after Equation~(2.4)).
   
   Pluging \eqref{eq:bbnew} into \eqref{eq:bborig} and taking the supremum over $\tf$ completes the 
   proof for general~$f$.
   If $f$ is hereditary, we have $f_{N_x^c \cup \{x\}}(\varrho + \delta_x) = g(x; \varrho) 
   f_{N_x^c}(\varrho)$, so that the additional representation of $\bphi$ claimed in \eqref{eq:bphiher}
   follows from the third line of Inequality~\eqref{eq:bbnew}.
\end{proof}

\section{The distance bounds} \label{sec:bound}

We begin this section by presenting the general setting for our main result, Theorem~\ref{thm:main}, partly compiling assumptions that were already mentioned, partly adding more specialized notation and conditions. Thereafter the main result and a number of corollaries are stated, and in the last subsection the corresponding proofs are given.

\subsection{Setting for Theorem~\ref{thm:main}} \label{ssec:setting}

Let $\xi$ be a point process on the compact metric space $(\mcx,\dzero)$ which has density $f$ with respect to $P_1$ and finite expectation measure $\mu = \mu_1$. Furthermore, let $\pi = \bigl( \pi(\cdot,x) ; \, x \in \mcx \bigr)$ be a $[0,1]$-valued random field. We assume that $\pi$ when interpreted as a random function on $\mcx$ takes values in a space $E \subset [0,1]^{\mcx}$ which is what we call an evaluable path space (i.e.\ $\Phi: E \times \mcx \to [0,1], (p,x) \mapsto p(x)$ is measurable) and that there exists a regular conditional distribution of $\pi$ given the value of $\xi$.
Neither of these assumptions presents a serious obstacle; we refer to Appendix~\ref{app:lepsrcd} for details. Let then $\Lambda$ be the
random measure given by $\Lambda(A) := \int_A \pi(x) \; \xi(dx)$ for $A \in \mcb$, and set $\Lambda^{[2]}(B) := \int_B \pi(x) \pi(\tx) \; \xi^{[2]} \bigl( d(x,\tx) \bigr)$ for any $B \in \mcb^2$.
 
Choose a neighborhood structure $(N_x)_{x \in \mcx}$ that is measurable in the sense of Condition~\eqref{eq:nhcondition}. We assume for every $x \in \mcx$ that $\pi(x)$ and $\pi \vert_{N_x^c}$ are both strictly locally dependent on $\xi$ in such a way that the corresponding ``regions of dependence'' do not interfere with one another. More exactly, this means the following: introduce an arbitrary metric $\tdzero$ on $\mcx$ that generates the same topology as $\dzero$,
and write $\BB(x,r)$ for the closed $\tdzero$-ball at $x \in \mcx$ with radius
$r \geq 0$ and $\BB(A,r) := \{ y \in \mcx ; \, \tdzero(y,x) \leq r \text{ for some $x \in A$} \}$ for the $\tdzero$-halo set of distance $r \geq 0$ around $A \subset \mcx$. Suppose then that we can fix a real number $R \geq 0$ such that
\begin{equation} \label{eq:minimalnbhdcond}
   \BB(x,R) \cap \BB(N_x^c,R) = \emptyset \quad \text{(i.e.\ $N_x \supset \BB(x,2R)$)}
\end{equation}
and
\begin{equation} \label{eq:locdepcond}
   \pi(x) \independent_{\xi \vert_{\BB(x,R)}} \xi \vert_{\BB(x,R)^c} \, \text{ and } \, \pi \vert_{N_x^c} 
   \independent_{\xi \vert_{\BB(N_x^c,R)}} \xi \vert_{\BB(N_x^c,R)^c}
\end{equation}
for every $x \in \mcx$, where $X \independent_{Z} Y$ denotes conditional independence of $X$ and $Y$ given $Z$.
If $Z$ is almost surely constant, this is just the (unconditional) independence of $X$ and $Y$; in particular we require $\pi(x) \independent \xi$ if $\alpha(\BB(x,R))=0$.
Set $\Aint := \Aint(x) := \BB(x,R)$ and $\Aext := \Aext(x) :=  \BB(N_x^c,R)$, where we usually suppress the location $x$ when it is clear from the context. 

We introduce two functions to control the dependences in $(\xi, \pi)$. The function $\bbeta: \mcx \to \Rplus$ is given by
\begin{equation} \label{eq:bbeta}
   \bbeta(x) := \int_{\mfn(\Aint)} \EE \bigl( \pi(x) \bigm| \xi \vert_{\Aint} = \varrhoint + \delta_x  
      \bigr) \int_{\mfn(\Aext)} \bigabs{ \fbar(\varrhoext, \varrhoint + \delta_x) } \;
      P_{1,\Aext}(d \varrhoext) \; P_{1,\Aint}(d \varrhoint), 
\end{equation}
where $\fbar(\varrhoext, \varrhoint + \delta_x) := f_{\Aext \cup \Aint}(\varrhoext + \varrhoint +   
\delta_x) - f_{\Aext}(\varrhoext) f_{\Aint}(\varrhoint + \delta_x)$, 
and hence controls the long range dependence within $\xi$, as well as the short range dependence of $\pi$
on~$\xi$. If $\alpha(\Aint) = 0$, the conditional expectation above is to be interpreted as $\EE( \pi(x) )$ for every $\varrhoint \in \mfn(\Aint)$. 
The function $\bgamma: \mcx \to \Rplus$ is taken to be a measurable function that satisfies
\begin{equation} \label{eq:bgamma}
   \int_{\mfn} \; \esssup_{C \in \sigma(\pi \vert_{N_x^c})} \Bigabs{ \cov \bigl( \pi(x), \indi{C} 
      \bigm| \xi = \varrho + \delta_x \bigr)} f(\varrho + \delta_x) \; P_1(d\varrho) \leq \bgamma(x),
\end{equation}
and hence controls the average long range dependence within $\pi$ given $\xi$. For the definition of the essential supremum of an arbitrary set of measurable functions (above, the functions $\bigl[ \varrho \mapsto \bigabs{\cov ( \pi(x), \indi{C} \mvert \xi = \varrho + \delta_x)} \bigr]$ for $C \in \sigma(\pi \vert_{N_x^c})$), see \cite{neveu65}, Proposition~II.4.1.

\subsection{Results}

We are now in the position to formulate our main theorem.
\begin{thm} \label{thm:main}
   Suppose that the assumptions of Subsection~\ref{ssec:setting} hold and
   write $\mu\obenpi = \EE \Lambda$ for the expectation measure of $\xi_{\pi}$.
   We then have
   \begin{equation*}
   \begin{split}
      \text{(i) } \, &\dtv \bigl( \msl(\xi_{\pi}), \poisson(\mu\obenpi) \bigr) \leq \bigl( \EE
         \Lambda \bigr)^2 \bigl( N(\mcx) \bigr) + \EE \Lambda^{[2]} \bigl( N(\mcx) \bigr) + \int_{\mcx}
         \bbeta(x) \; \alpha(dx) + 2 \int_{\mcx} \bgamma(x) \; \alpha(dx); \\[2mm]
      \text{(ii) } \, &\dtwo \bigl( \msl(\xi_{\pi}), \poisson(\mu\obenpi) \bigr) \leq M_2 \bigl(
         \mu\obenpi \bigr) \Bigl( \bigl( \EE \Lambda \bigr)^2 \bigl( N(\mcx) \bigr) + \EE \Lambda^{[2]} 
         \bigl( N(\mcx) \bigr) \Bigr) \\
      &\hspace*{45.5mm} {} + M_1 \bigl( \mu\obenpi \bigr) \biggl( \int_{\mcx}
         \bbeta(x) \; \alpha(dx) + 2 \int_{\mcx} \bgamma(x) \; \alpha(dx) \biggr);
   \end{split}
   \end{equation*}
   where
   \begin{equation*}
   \begin{split}
      M_1(\mu\obenpi) &= \min \biggl( 1, \frac{1.647}{\sqrt{\EE \Lambda(\mcx)}} \biggr), \text{ and} \\
      M_2(\mu\obenpi) &= \min \biggl( 1, \frac{11}{6 \, \EE \Lambda(\mcx)} \Bigl(1 + 2 \log^{+} \bigl( 
         \tfrac{6}{11} \EE \Lambda(\mcx) \bigr) \Bigr) \biggr).
   \end{split}
   \end{equation*}
\end{thm}
If $\xi$ and $\pi$ are independent, we obtain an interesting special case, where the upper bound can be expressed in terms of essentially the quantities appearing in Theorem~\ref{thm:bb}, which are based 
solely on $\xi$, and some rather straightforward quantities based on $\pi$.
\begin{cor}[Independent case]  \label{cor:main}
   Suppose that $\xi$ is a point process on $(\mcx, \dzero)$ which has density $f$ with respect to $P_1$ and
   finite expectation measure $\mu = \mu_1$. Let $\pi = \bigl( \pi(\cdot, x); \, x \in \mcx \bigr)$ be a
   $[0,1]$-valued random field that has an evaluable path space $E \subset [0,1]^{\mcx}$ and is independent
   of~$\xi$. 
   Choose an arbitrary neighborhood structure $(N_x)_{x \in \mcx}$ (always $x \in N_x$!) that is measurable 
   in the sense of Condition \eqref{eq:nhcondition}, and take $\bgamma: \mcx \to \Rplus$ to be a measurable
   function that satisfies
   \begin{equation} \label{eq:bgammai}
      \sup_{C \in \sigma(\pi \vert_{N_x^c})} \bigabs{ \cov( \pi(x), \indi{C} )} \cdot u(x) \leq \bgamma(x)
   \end{equation}
   for almost every $x$.
   Note that the expectation measure of $\xi_{\pi}$ is $\mu\obenpi(\cdot) = \EE \Lambda(\cdot) =
   \int_{\bs{\cdot}} \EE \pi(x) \; \mu_1(dx)$, and let $M_1(\mu\obenpi)$ and $M_2(\mu\obenpi)$ be defined
   as in Theorem~\ref{thm:main}. 
   We then obtain
   \begin{equation*}
   \begin{split}
      \text{(i) } \, &\dtv \bigl( \msl(\xi_{\pi}), \poisson(\mu\obenpi) \bigr) \\[1mm]
      &\hspace*{10mm} \leq \int_{\mcx} 
         \int_{N_x} \EE \pi(x) \EE \pi(\tx) \; \mu_1(d\tx) \; \mu_1(dx) +
         \int_{N(\mcx)} \EE \bigl( \pi(x) \pi(\tx) \bigr) \; \mu_{[2]} \bigl( d(x,\tx) \bigr) \\
      &\hspace*{18mm} {} + \int_{\mcx} \EE \pi(x) \, \bphi(x)\; \alpha(dx) + 2 \int_{\mcx} \bgamma(x) \;
         \alpha(dx); \\[2mm]
      \text{(ii) } \, &\dtwo \bigl( \msl(\xi_{\pi}), \poisson(\mu\obenpi) \bigr) \\[1mm]
      &\leq M_2 \bigl( \mu\obenpi \bigr) \biggl(
         \int_{\mcx} \int_{N_x} \EE \pi(x) \EE \pi(\tx) \; \mu_1(d\tx) \; \mu_1(dx) +
         \int_{N(\mcx)} \EE \bigl( \pi(x) \pi(\tx) \bigr) \; \mu_{[2]} \bigl( d(x,\tx) \bigr) \biggr) \\
      &\hspace*{7mm} {} + M_1 \bigl( \mu\obenpi \bigr) \biggl( \int_{\mcx} \EE \pi(x) \, \bphi(x) \;
         \alpha(dx) + 2 \int_{\mcx} \bgamma(x) \; \alpha(dx) \biggr); 
   \end{split}
   \end{equation*}
   where $\bphi$ is given by Equation~\eqref{eq:bphinonher} and, if $f$ is hereditary, by Equation~\eqref
   {eq:bphiher}.
   \end{cor}
A further corollary is given for the case of a deterministic retention field, which means that the retention decisions are independent of each other given the point process $\xi$. We formulate only a very special case,
where the point process lives on $\RD$ and various spatial homogeneities are assumed, which leads to a particularly simple upper bound.
The corollary also illustrates how we can deal with the issue of boundary effects in the neighborhood structure by extending the $N_x$ beyond $\mcx$ in such a way that they are translated versions of each other and that the same holds true for their ``complements'' $M_x \setminus N_x$.

In this article we always tacitly assume that $\RD$ is equipped with the Euclidean topology.
Write $\abs{A} := \LebD(A)$ for any Borel set $A \subset \RD$, and $\mcx + \mcx' := \{x + x'; \, x \in \mcx, x' \in \mcx' \}$ and $\mcx - x := \{\tx - x; \, \tx \in \mcx \}$ for $\mcx, \mcx' \subset \RD$ and $x \in \RD$. For the definition of point processes on non-compact spaces and elementary concepts such as stationarity, we refer to the standard point process literature (e.g.\ \cite{kallenberg86} or \cite{dvj88}).
\begin{cor}[Constant case] \label{cor:mainconst}
   Let $\mcx, \mcy \subset \RD$ be two compact sets, where $\mcx$ has positive volume and $\mcx + \bigcup_{x 
   \in \mcx} (\mcx - x) \subset \mcy$, and consider a metric $\dzero \leq 1$ on $\mcx$ that generates 
   the Euclidean topology (typically, $\dzero(x,y) = \min(\abs{x-y},1)$ for all $x,y \in \mcx$).

   \nin
   Suppose that $\zeta$ is a stationary point process on $\RD$ whose restriction $\xi := \zeta \vert_{\mcy}$
   is a Gibbs process with density $f$ with respect to $P_1  = \poisson(\LebD \vert_{\mcy})$ and
   finite expectation measure $\mu = \mu_1 = m_1 \LebD \vert_{\mcy}$. Denote by $\mck$ the second reduced
   moment measure of $\zeta$ (see around Equation~\eqref{eq:kmeasure} for details).  
   Let $(N_x)_{x \in \mcx}$ and $(M_x)_{x \in \mcx}$ 
   be two neighborhood structures whose sets $N_x := x + N$ and $M_x := x + M$ (not necessarily $\subset \mcx$ 
   now!) are translated copies of single bounded measurable sets $N, M \subset \RD$ which are chosen
   in such a way that $N \subset M$ and $\mcx \subset M_x \subset \mcy$ for every $x \in \mcx$.
   Choosing our retention field $\pi \equiv p \in [0,1]$ to be deterministic and constant and noting that
   $\mu\obenp \vert_{\mcx} = p \, m_1 \LebD \vert_{\mcx}$,
   we then have
   \begin{equation*}
   \begin{split}
      \text{(i) } \, &\dtv \bigl( \msl(\xi_{p} \vert_{\mcx}), \poisson(p \, m_1 \LebD \vert_{\mcx}) \bigr)
         \leq p^2
         m_1^2 \abs{\mcx} (\abs{N} + \mck(N))+ p \abs{\mcx} \, \EE \abs{\Gamma - \EE \Gamma}; \\[2mm]
      \text{(ii) } \, &\dtwo \bigl( \msl(\xi_{p} \vert_{\mcx}), \poisson(p \, m_1 \LebD \vert_{\mcx}) \bigr) 
         \leq \min \Bigl( p \, m_1 \abs{\mcx}, \tfrac{11}{6} \bigl(1 + 2 \log^{+} \bigl( \tfrac{6 p \, m_1
         \abs{\mcx}}{11}
         \bigr) \bigr) \Bigr) p \, m_1 (\abs{N} + \mck(N)) \\
      &\hspace*{55mm} {} + \min \bigl( \sqrt{p \, m_1 \abs{\mcx}}, 1.647 \bigr) \sqrt{p \, m_1 \abs{\mcx}} 
         \tfrac{1}{m_1} \EE \abs{\Gamma - \EE \Gamma}, 
   \end{split}
   \end{equation*}
   where $\Gamma$ is an arbitrary random variable that has the same distribution as $g'(x; \xi \vert_{M_x 
   \setminus N_x}) = f_{(M_x \setminus N_x) \cup \{x\}}(\xi \vert_{M_x \setminus N_x} + \delta_x) \big/ f_{M_x 
   \setminus N_x}(\xi \vert_{M_x \setminus N_x})$ (for one and therefore every $x \in \mcx$). 
\end{cor}
\begin{rem}  \label{rem:mainconst}
   While it seems very appealing to admit $M_x = \mcy = \RD$, this case actually requires a different and 
   somewhat technically involved construction for the conditional density $g'(x; \xi \vert_{\RD
   \setminus N_x})$, because it cannot reasonably be assumed that a density of a point process distribution 
   with respect to the standard Poisson process distribution exists if the state space is $\RD$ (consider
   for example a hard core process: the hard core event that no two 
   points are closer than some fixed distance $r>0$ is a $\poisson(\LebD)$-null set).
   As a matter of fact, the natural setting is that of a Gibbs process on the whole of $\RD$ defined via a 
   stable potential on the set of finite point measures on $\RD$, which essentially
   provides us with ``compatible'' conditional densities for the point process distribution on bounded Borel
   sets given the point process outside (see~\cite{preston76}, from page 6.9 onwards, for a detailed 
   construction).
   For a fixed bounded Borel 
   set $\mcy \subset \RD$ we write $f_B(\cdot \mvert \tau) : \mfn(B) \to \Rplus$ for the conditional density of 
   $\xi \vert_B$ given $\xi \vert_{\mcy^c} = \tau$. It can then be seen that the crucial inequality
   \begin{equation*}
      \EE \bigabs{\tg(x; \xi \vert_{\mcx \setminus N_x}) - m_1} \leq \EE \bigabs{g'(x; \xi \vert_{M_x \setminus 
       N_x}) - m_1} 
   \end{equation*}
   (see Inequality~\eqref{eq:extendingg}) and hence the proof of Corollary~\ref{cor:mainconst} can be 
   reproduced under very general conditions if $M_x = \RD$,
   where
   \begin{equation*}
      \tg(x; \tvarrho) = \frac{\int_{\mfn(\mcy^c)} f_{(\mcx \setminus N_x) \cup \{x\}}(\tvarrho + \delta_x 
      \mvert \tau) \; \PP (\xi \vert_{\mcy^c})^{-1}(d \tau)}{\int_{\mfn(\mcy^c)} f_{\mcx \setminus
      N_x}(\tvarrho \mvert \tau) \; \PP (\xi \vert_{\mcy^c})^{-1}(d \tau)}
   \end{equation*}
   for every $\tvarrho \in \mfn(\mcx \setminus N_x)$, and
   \begin{equation*}
      g'(x; \varrho) = \frac{f_{(\mcy \setminus N_x) \cup \{x\}}(\varrho \vert_{\mcy \setminus N_x} + \delta_x 
      \mvert \varrho \vert_{\mcy^c})}{f_{\mcy \setminus N_x}(\varrho \vert_{\mcy \setminus N_x} 
      \mvert \varrho \vert_{\mcy^c})}
   \end{equation*}
   for every $\varrho \in \mfn(\RD \setminus N_x)$ (as earlier we set such fractions to zero if the denominator
   is zero).
   By the construction of the Gibbs process on $\RD$ (using Equation~(6.10) in \cite{preston76}), the term
   $g'(x; \varrho)$ does not depend on the choice of $\mcy \supset N_x$, except for $\varrho$ in a $\PP(\xi 
   \vert_{\RD \setminus N_x})^{-1}$-null set. It can be interpreted as the conditional density of a point at
   $x$ given $\xi \vert_{\RD \setminus N_x} = \varrho$.
\end{rem}
In the next result, the situation of Theorem~\ref{thm:main} and its corollaries is examined for the case where we compensate for the thinning by contracting the state space as it was done in \cite{schumi2}.
\begin{cor}[Thinning and contraction in $\bs{\RD}$] \label{cor:maincontracted}
   Suppose that $\mcx$ is a compact subset of $\RR^D$ and that $\alpha = \LebD \vert_{\mcx}$. Let $T \geq 1$
   and $\kappa_T: \RD \to \RD, x \mapsto \frac{1}{T}x$. Assume furthermore
   that the metric $\dzero$ on $\mcx$ generates the Euclidean topology and 
   satisfies $\dzero \bigl( \kappa_T(x), \kappa_T(y) \bigr) \leq \dzero(x,y)$ for every $x,y \in \mcx$.\\
   Then Theorem~\ref{thm:main} and Corollaries~\ref{cor:main} and~\ref{cor:mainconst} remain true if 
   the point processes on the left hand sides of the estimates are replaced by their respective 
   image processes under the contraction~$\kappa_T$. We thus obtain under the general prerequisites of
   Theorem~\ref{thm:main}
   \begin{equation*}
   \begin{split}
      \text{(i) } \, &\dtv \bigl( \msl(\xi_{\pi} \kappa_T^{-1}), \poisson(\mu\obenpi \kappa_T^{-1}) \bigr) 
         \leq \bigl( \EE \Lambda \bigr)^2 \bigl( N(\mcx) \bigr) + \EE \Lambda^{[2]} \bigl( 
         N(\mcx) \bigr) + \int_{\mcx} \bbeta(x) \; dx + 2 \int_{\mcx} \bgamma(x) \; dx; 
         \\[2mm]
      \text{(ii) } \, &\dtwo \bigl( \msl(\xi_{\pi} \kappatinv), \poisson(\mu\obenpi \kappatinv) \bigr)
         \leq M_2 \bigl( \mu\obenpi \bigr) \Bigl( \bigl( \EE \Lambda \bigr)^2 \bigl( N(\mcx) \bigr)
         + \EE \Lambda^{[2]} \bigl( N(\mcx) \bigr) \Bigr) \\
      &\hspace*{58.5mm} {} + M_1 \bigl( \mu\obenpi \bigr) \biggl( \int_{\mcx}
         \bbeta(x) \; dx + 2 \int_{\mcx} \bgamma(x) \; dx \biggr).
   \end{split}
   \end{equation*}
\end{cor}
\begin{rem}[Comparison with \cite{schumi2}]  \label{rem:compmain}
   The setting of Corollary~\ref{cor:maincontracted} corresponds in large parts to the situation
   in \cite{schumi2}, especially if we set $\mcx := \kappatinv(J)$ for a fixed compact set $J \subset \RD$   
   and compare our statement~(ii) to Theorem~3.B.
   
   \nin
   It is more strict in essentially two respects. First, of course, 
   we admit only point processes whose distributions are absolutely continuous with respect to a homogeneous
   Poisson process.  
   Secondly, we require strict local dependence of $\pi$ on $\xi$ (see Condition~\eqref{eq:locdepcond}),
   which in \cite{schumi2} was only done for Section~4 (in slightly different form), but which also 
   gives the direct benefit of a conceptually simpler and more intuitive control of the long range 
   dependences.
   
   \nin
   On the other hand, the setting of Corollary~\ref{cor:maincontracted} gives us more freedom 
   than we had in \cite{schumi2} in the sense that the objects live on a general compact subset of 
   $\RD$, that there are only minimal conditions on the moment measures (as opposed to Assumption~1 in
   \cite{schumi2}), and that the choice of $\dzero$ and of the neighborhoods of
   strong dependence $N_x$ is much wider.
   
   \nin
   Regarding the upper bounds obtained we have clearly improved. The terms in our statement~(ii) above have 
   their counterparts in the various terms in Theorem~3.B of \cite{schumi2} (with the integrals over
   $\bbeta(x)$ and $\bgamma(x)$ being summerized as a single long range dependence term), but have become
   simpler and more natural, 
   without any suprema or infima over discretization cuboids and with explicit and manageable constants.
   The bound as a whole is somewhat better (if we use a heuristic approach for comparing the long range 
   dependence terms) and quite a bit more easily applied, which can be seen from comparing the application in
   Subsection~\ref{ssec:ibooleanapplic} with the same application in Subsection~3.3 of \cite{schumi2}.
\end{rem}
\begin{rem}
   As pointed out in \cite{schumi2}, it may be desirable to approximate the distribution of the thinned
   point process by a Poisson process distribution that has a somewhat different expectation
   measure. Corresponding distance estimates can easily be obtained from upper bounds
   for distances between Poisson process distributions. We have $\dtv \bigl( \poisson(\lambda), \poisson(\mu) 
   \bigr) \leq \norm{\lambda - \mu}$ for $\lambda, \mu \in \mfm$ by Remark~2.9 in \cite{bb92}, where
   $\norm{\cdot}$ denotes the total variation norm for signed measures. For $\dtwo$ an
   upper bound is given as Inequality~(A.3) in \cite{schumi2} (which is the same as Inequality~(2.8)
   in \cite{brownxia95}).
\end{rem}

\subsection{Proofs} \label{ssec:proofs}

\begin{proof}[Proof of Theorem~\ref{thm:main}]

By Lemma~\ref{lem:fp} in the appendix a density $f\obenpi$ of $\xi_{\pi}$ with respect to $P_1$ exists, and the finiteness of $\mu$ implies the finiteness of $\mu\obenpi$. Hence we can apply Theorem~\ref{thm:bb}.

The integrals in the upper bound can be further evaluated as follows. For the first two integrals (basic term and strong dependence term), we have by Lemma~\ref{lem:mms} that
\begin{equation} \label{eq:bttotal}
   \int_{\mcx} \mu\obenpi(N_x) \; \mu\obenpi(dx)
   = \int_{\mcx} \bigl( \EE \Lambda \bigr)(N_x) \; \bigl( \EE \Lambda \bigr) (dx)
   = \bigl( \EE \Lambda \bigr)^2 \bigl( N(\mcx) \bigr)
\end{equation}
and
\begin{equation} \label{eq:sdttotal}
   \EE \biggl( \int_{\mcx} \xi_{\pi}(\dot{N}_x) \; \xi_{\pi}(dx) \biggr) = \EE \bigl( \xi_{\pi}^{[2]} 
   \bigl( N(\mcx) \bigr) \bigr) = \EE \bigl( \Lambda^{[2]} \bigl( N(\mcx) \bigr) \bigr). 
\end{equation}

For the third integral (weak dependence term) some more work is necessary. The term that we would like to estimate is
\begin{equation} \label{eq:wdtgoal}
\begin{split}
   2 \int_{\mcx} \sup_{C \in \mcn(N_x^c)} \biggabs{\int_C \bigl[ f\obenpi_{N_x^c \cup \{x\}}(\varrho + 
   \delta_x) - f\obenpi_{N_x^c}(\varrho) u\obenpi(x) \bigr] \; P_{1,N_x^c}(d \varrho)} \; \alpha(dx),
\end{split}
\end{equation}
where $u\obenpi$ is the density of $\mu\obenpi$.
Equations~\eqref{eq:wdt_crucial0} and~\eqref{eq:wdt_crucial1} from the appendix imply that, for almost every
$x \in \mcx$ and for $C \in \mcn(N_x^c)$,
\begin{equation} \label{eq:wdtstart0}
      \int_C f_{N_x^c}\obenpi(\varrho) \; P_{1,N_x^c}(d \varrho) = \int_{\mfn} \EE \bigl(
      Q_{C}\obenpi(\sigma \vert_{N_x^c}) \bigm| \xi = \sigma \bigr) f(\sigma) \; P_1(d \sigma)
\end{equation}
and
\begin{equation} \label{eq:wdtstart1}
   \int_C f_{N_x^c \cup \{x\}}\obenpi(\varrho + \delta_x) \; P_{1,N_x^c}(d \varrho) = \int_{\mfn}
      \EE \bigl( \pi(x) Q_{C}\obenpi(\sigma \vert_{N_x^c})
      \bigm| \xi = \sigma + \delta_x \bigr) f(\sigma + \delta_x) \; P_1(d \sigma),
\end{equation}
where $Q_C\obenp(\sigma) = \sum_{\varrho \subset \sigma, \varrho \in C} \bigl( \prod_{s \in \varrho} p(s) \bigr) \bigl( \prod_{\ts \in \sigma \setminus \varrho} (1-p(\ts)) \bigr)$ for every $\sigma \in \mfn^{*}$ and every $p \in E$, so that $\bigl[(p, \sigma) \mapsto Q_C\obenp(\sigma)\bigr]$ is $\mce \otimes \mcn$-measurable.
By Equation~\eqref{eq:pregnzconsequence} we have furthermore that
\begin{equation}
   u\obenpi(x) = \int_{\mfn(N_x^c)} f_{N_x^c \cup \{x\}}\obenpi
   (\varrho + \delta_x) \; P_{1, N_x^c}(d \varrho) = \int_{\mfn} \EE \bigl( \pi(x) \bigm| \xi = \sigma
   + \delta_x \bigr) f(\sigma + \delta_x) \; P_1(d \sigma),
\end{equation}
using that $Q_{\mfn(N_x^c)}\obenp(\sigma) = 1$ for every $\sigma \in \mfn^{*}(N_x^c)$. 

The absolute value term in~\eqref{eq:wdtgoal} can then be estimated as
\label{page:monster}
\begin{equation} \label{eq:monster}
\begin{split}
   \biggabs{ &\int_C \bigl[ f_{N_x^c \cup \{x\}}\obenpi(\varrho + \delta_x) - f_{N_x^c}\obenpi (\varrho)
      u\obenpi(x) \bigr] \; P_{1, N_x^c}(d \varrho) } \\[3mm]
   &= \biggabs{\int_{\mfn} \EE \bigl( \pi(x) Q_C\obenpi(\sigma \vert_{N_x^c}) \bigm| \xi
      = \sigma + \delta_x \bigr) f(\sigma + \delta_x) \; P_{1}(d \sigma) \\
   &\hspace*{5mm} {} -
      \int_{\mfn} \EE \bigl( Q_C\obenpi(\sigma \vert_{N_x^c}) \bigm| \xi = \sigma \bigr)
      f(\sigma) \; P_{1}(d \sigma) \, \cdot \, \int_{\mfn} \EE \bigl( \pi(x) 
      \bigm| \xi = \sigma + \delta_x \bigr) f(\sigma +
      \delta_x) \; P_{1}(d \sigma) } \\[3mm]
   &\leq \biggabs{ \int_{\mfn} \EE \bigl( \pi(x) Q_C\obenpi(\sigma \vert_{N_x^c}) \bigm| \xi
      = \sigma + \delta_x \bigr) f(\sigma + \delta_x) \; P_{1}(d \sigma) \\
   &\hspace*{10mm} {} - \int_{\mfn} \EE \bigl( Q_C\obenpi(\sigma \vert_{N_x^c}) \bigm| \xi = \sigma +
      \delta_x \bigr) \: \EE \bigl( \pi(x) \bigm| \xi = \sigma + \delta_x \bigr) f(\sigma +
      \delta_x) \; P_{1}(d \sigma) } \\[1mm]
   &\hspace*{5mm} {} + \biggabs{ \int_{\mfn} \EE \bigl( Q_C\obenpi(\sigma \vert_{N_x^c}) \bigm| \xi = \sigma
      + \delta_x \bigr) \: \EE \bigl( \pi(x) \bigm| \xi = \sigma + \delta_x \bigr) f(\sigma +
      \delta_x) \; P_{1}(d \sigma) \\
   &\hspace*{10mm} {} - \int_{\mfn} \EE \bigl( Q_C\obenpi(\sigma \vert_{N_x^c}) \bigm| \xi = \sigma \bigr)
      f(\sigma) \; P_{1}(d \sigma) \, \cdot \, \int_{\mfn} \EE \bigl( \pi(x) 
      \bigm| \xi = \sigma + \delta_x \bigr) f(\sigma +
      \delta_x) \; P_{1}(d \sigma) } \\[3mm]
   &= \biggabs{ \int_{\mfn} \cov \bigl( \pi(x), Q_C\obenpi \bigl( \sigma \vert_{N_x^c} \bigr) \bigm|
      \xi = \sigma + \delta_x \bigr) f(\sigma + \delta_x) \; P_1(d \sigma) } \\[1mm]
   &\hspace*{10mm} {} + \biggabs{ \int_{\mfn(\Aext)} \int_{\mfn(\Aint)} \EE \bigl( Q_C\obenpi(\varrhoext
      \vert_{N_x^c}) \bigm| \xi \vert_{\Aext}= \varrhoext
      \bigr) \: \EE \bigl( \pi(x) \bigm| \xi \vert_{\Aint} = \varrhoint + \delta_x \bigr) \\[-1.5mm]
   &\hspace*{60mm} f_{\Aext \cup \Aint}(\varrhoext +
      \varrhoint + \delta_x) \; P_{1, \Aint} (d \varrhoint) \; P_{1, \Aext} (d \varrhoext) \\[1mm]
   &\hspace*{16mm} {} - \int_{\mfn(\Aext)} \EE \bigl( Q_C\obenpi(\varrhoext \vert_{N_x^c}) \bigm|
      \xi \vert_{\Aext}= \varrhoext \bigr)
      f_{\Aext}(\varrhoext) \; P_{1, \Aext}(d \varrhoext) \\
   &\hspace*{21mm} {} \cdot \int_{\mfn(\Aint)} \EE
      \bigl( \pi(x) \bigm| \xi \vert_{\Aint} = \varrhoint + \delta_x \bigr) f_{\Aint}(\varrhoint +
      \delta_x) \; P_{1, \Aint}(d \varrhoint) }, \\[-13.5mm]
   &\hspace*{2mm}
\end{split}
\end{equation}
where Condition~\eqref{eq:locdepcond} was used for the last equality. Note that $Q_C\obenpi(\sigma \vert_{N_x^c})$ depends on $\pi$ only via $\pi \vert_{N_x^c}$, and does so in a $\mce(N_x^c)$-measurable way, where $\mce(N_x^c) = \sigma \bigl( [ \tp \mapsto \tp(x)]; \, x \in N_x^c \bigr)$ is the canonical $\sigma$-algebra on $E(N_x^c) := \{ p \vert_{N_x^c}; \, p \in E \}$.

The first summand on the right hand side of Inequality~\eqref{eq:monster} can then be bounded further as
\begin{equation*}
\begin{split}
   \biggabs{ \int_{\mfn} \cov \bigl( \pi(x), Q_C\obenpi &\bigl( \sigma \vert_{N_x^c} \bigr) \bigm|
      \xi = \sigma + \delta_x \bigr) f(\sigma + \delta_x) \; P_1(d \sigma) } \\
   &\leq \int_{\mfn} \esssup_{\hspace*{0.5pt} h: E(N_x^c) \to [0,1]} \Bigabs{ \cov \bigl( \pi(x), h(\pi
      \vert_{N_x^c} \bigr) \bigm| \xi = \sigma + \delta_x \bigr) } f(\sigma + \delta_x) \; P_1(d \sigma) \\
   &= \int_{\mfn} \esssup_{\hspace*{0.5pt} h: E(N_x^c) \to \{0,1\}} \Bigabs{ \cov \bigl( \pi(x), h(\pi
      \vert_{N_x^c} \bigr) \bigm| \xi = \sigma + \delta_x \bigr) } f(\sigma + \delta_x) \; P_1(d \sigma) \\
   &\leq \bgamma(x),
\end{split}
\end{equation*}
where the essential suprema are taken over all $\mce(N_x^c)$-measurable functions with values in $[0,1]$ and $\{0,1\}$, respectively. The third line is obtained by
\begin{equation*}
\begin{split}
   \bigabs{ \cov(X,Y) } &\leq \linftynorm{Y} \bigabs{ \cov \bigl( X, \sign(\tY) \bigr) } \\
   &\leq 2 \linftynorm{Y} \max \Bigl( \bigabs{\cov \bigl( X, {\rm I}[\tY > 0] \bigr)}, \bigabs{\cov \bigl( X, 
      {\rm I}[\tY < 0] \bigr)} \Bigr)
\end{split}
\end{equation*}
for all random variables $X \in L_1$ and $Y \in L_{\infty}$, and for $\tY := \EE(X \mvert Y) - \EE  X$ (see \cite{doukhan94}, Section~1.2, proof of Lemma~3), where we set $X := \pi(x)$ and $Y := h(\pi \vert_{N_x^c}) - 1/2$.

For the second summand on the right hand side of Inequality~\eqref{eq:monster}, we use the notation $F_C(\varrhoext) := \EE \bigl( Q_C\obenpi(\varrhoext \vert_{N_x^c}) \bigm| \xi \vert_{\Aext}= \varrhoext \bigr)$ and $G(\varrhoint + \delta_x) := \EE \bigl( \pi(x) \bigm| \xi \vert_{\Aint} = \varrhoint + \delta_x \bigr)$, and bound it as
\begin{equation*}
\begin{split}
    &\biggabs{ \int_{\mfn(\Aext)} \int_{\mfn(\Aint)}
      F_C(\varrhoext) G(\varrhoint + \delta_x) \fbar(\varrhoext, \varrhoint + \delta_x)
      \; P_{1, \Aint} (d \varrhoint) \; P_{1, \Aext}(d \varrhoext)} \\
   &\hspace*{1mm} \leq \int_{\mfn(\Aint)} G(\varrhoint + \delta_x) \sup_{F: \mfn(\Aext) \to 
      [0,1]} \biggabs{
      \int_{\mfn(\Aext)} F(\varrhoext) \fbar(\varrhoext, \varrhoint + \delta_x) \; P_{1, \Aext}(d 
      \varrhoext) } \; P_{1, \Aint} (d \varrhoint) \\[1mm]
   &\hspace*{1mm}= \frac12 \int_{\mfn(\Aint)} G(\varrhoint + \delta_x) \int_{\mfn(\Aext)} \bigabs{ \fbar
      (\varrhoext, \varrhoint + \delta_x) } \; P_{1, \Aext}(d \varrhoext) \; P_{1, \Aint} (d \varrhoint),
\end{split}
\end{equation*}
where the supremum is taken over $\mcn(\Aext)$-measurable functions. 
The third line is shown by setting $F_0(\varrhoext) := {\rm I} \bigl[ \fbar(\varrhoext,\varrhoint+\delta_x) > 0 \bigr]$ and noting that $F_0: \mfn(\Aext) \to [0,1]$ is measurable and maximizes the absolute value
term in the second line.

Thus the total estimate for the weak dependence term is
\begin{equation} \label{eq:wdttotal}
\begin{split}
   2 \int_{\mcx} \sup_{C \in \mcn(N_x^c)} \biggabs{\int_C \bigl[ &f\obenpi_{N_x^c \cup \{x\}}(\varrho + 
      \delta_x) - f\obenpi_{N_x^c}(\varrho) u\obenpi(x) \bigr] \; P_{1,N_x^c}(d \varrho)} \; \alpha(dx) \\
   &\hspace*{30mm} \leq 2 \int_{\mcx} \bgamma(x) \; \alpha(dx) + \int_{\mcx} \bbeta(x) \; \alpha(dx).
\end{split} 
\end{equation}
Pluging \eqref{eq:bttotal}, \eqref{eq:sdttotal}, and \eqref{eq:wdttotal} into Theorem~\ref{thm:bb} yields statement~(i), and, since $\abs{\mu\obenpi} = \EE \Lambda(\mcx)$, also statement~(ii). 
\end{proof}
\begin{proof}[Proof of Corollary~\ref{cor:main}]
   We aim at applying Theorem~\ref{thm:main} for $R=0$. Clearly Condition~\eqref{eq:minimalnbhdcond} holds for 
   any neighborhood structure. By the independence of $\xi$ and $\pi$ we have $\msl(\pi)$ as a 
   regular conditional distribution of $\pi$ given the value of $\xi$ and we see that
   Condition~\eqref{eq:locdepcond} is satisfied, that $\bbeta(x) = \EE(\pi(x)) \bphi(x)$ for almost every $x$, 
   and that Inequality~\eqref{eq:bgamma} simplifies to \eqref{eq:bgammai} by Equation~\eqref{eq:u}.
   Using the representation $\xi = \sum_{i=1}^V \delta_{S_i}$ from the proof of
   Lemma~\ref{lem:mms}, we have furthermore that
   \begin{equation*}
   \begin{split}
      \EE \Lambda(A) &= \EE \biggl( \EE \biggl( \sum_{i=1}^V \pi(S_i) \hspace*{1pt} {\rm I}[S_i \in A] 
         \biggm| \xi \biggr) \biggr) \\
      &= \EE \biggl( \sum_{i=1}^V \EE \bigl( \pi(S_i) \bigm| \xi \bigr) \hspace*{1pt} {\rm I}[S_i \in A] 
         \biggr) \\
      &= \EE \biggl( \sum_{i=1}^V \bigr( \EE \pi(x) \bigr)\big\vert_{x=S_i} \hspace*{0.5pt} {\rm I}[S_i \in A] 
         \biggr) \\
      &= \int_A \EE \pi(x) \; \mu_1(dx)
   \end{split}
   \end{equation*}
   for every $A \in \mcb$, and by the analogous computations that 
   \begin{equation*}
      \EE \Lambda^{[2]}(B) = \int_{B} \EE \bigl(\pi(x) \pi(\tx) \bigr) \; \mu_{[2]} \bigl( d(x,\tx) \bigr)
   \end{equation*}
   for every $B \in \mcb^2$. Based on these results we can apply Theorem~\ref{thm:main} and obtain the upper 
   bounds stated.
\end{proof}
\begin{proof}[Proof of Corollary~\ref{cor:mainconst}]
We apply Corollary~\ref{cor:main} for the point process $\txi := \xi \vert_{\mcx}$, which has hereditary density
$\tf := f_{\mcx}$ with respect to $\poisson(\LebD \vert_{\mcx})$ and expectation measure $\tmu = m_1 \LebD \vert_{\mcx}$, where all of these objects are interpreted as living on $\mcx$ (as opposed to living on $\mcy$ and being trivial on $\mcy \setminus \mcx$). Consider as neighborhood structure $(\tN_x)_{x \in \mcx}$ given by $\tN_x := N_x \cap \mcx$, write $\tN_x^c$ for the complement of $\tN_x$ in $\mcx$, and set $N(\mcx) := \{ (x,y) \in \mcx \times \mcy; \, y \in N_x \}$ and $\tN(\mcx) := \{ (x,y) \in \mcx^2; \, y \in \tN_x \}$, which are measurable by the fact that the $N_x$ are translated copies of a single measurable set.
Denoting the conditional density based on $\tf$ by $\tg$, we obtain for the $\bphi(x)$-term
\begin{equation} \label{eq:extendingg}
\begin{split}
   \EE \bigabs{&\tg(x; \txi \vert_{\tN_x^c}) - m_1} \\
   &\: = \int_{\mfn(\tN_x^c)} \bigabs{ \tf_{\tN_x^c \cup \{x\}}(\varrho + \delta_x) - \tf_{\tN_x^c}(\varrho)
      \, m_1 } \; P_{1,\tN_x^c}(d \varrho) \\
   &\: = \int_{\mfn(\mcx \setminus N_x)} \bigabs{ f_{(\mcx \setminus N_x) \cup \{x\}}(\varrho + \delta_x) -
      f_{\mcx \setminus N_x}(\varrho) \, m_1 } \; P_{1, \mcx \setminus N_x}(d \varrho) \\
   &\: = \int_{\mfn(\mcx \setminus N_x)} \biggabs{ \int_{\mfn(M_x \setminus (\mcx \cup N_x))} \bigl( f_{(M_x 
      \setminus N_x) \cup \{x\}}(\varrho + \tvarrho + \delta_x) \\[-3.5mm]
   &\: \hspace*{58mm} {} - f_{(M_x \setminus N_x)}(\varrho + \tvarrho)
      \, m_1 \bigr) \; P_{1, M_x \setminus (\mcx \cup N_x)}(d \tvarrho) } \; P_{1, \mcx \setminus N_x}(d 
      \varrho) \\ 
   &\: \leq \int_{\mfn(M_x \setminus N_x)} \bigabs{ f_{(M_x \setminus N_x) \cup \{x\}}(\sigma + \delta_x) -
      f_{M_x \setminus N_x}(\sigma) \, m_1 } \; P_{1, M_x \setminus N_x}(d \sigma) \\[1.5mm]
   &\: = \EE \bigabs{g'(x; \xi \vert_{M_x \setminus N_x}) - m_1}, 
\end{split}
\end{equation}
and thus by Corollary~\ref{cor:main} that
\begin{equation*}
\begin{split}
   \dtv \bigl( \msl( &\xi_p \vert_{\mcx}), \poisson(p \, m_1 \LebD \vert_{\mcx}) \bigr) \\[1.5mm]
   &= \dtv \bigl( \msl(\txi_p), \poisson(\tmu\obenp) \bigr) \\[0.5mm]
   &\leq p^2 \int_{\mcx} \mu_1(\tN_x) \; \mu_1(dx) + p^2 \mu_{[2]}(\tN(\mcx)) + p \int_{\mcx}
      \EE \bigabs{\tg(x; \txi \vert_{\tN_x^c}) - m_1} \; dx + 0 \\
   &\leq p^2 \int_{\mcx} \mu_1(N_x) \; \mu_1(dx) + p^2 \mu_{[2]}(N(\mcx)) + p \int_{\mcx}
      \EE \bigabs{g'(x; \xi \vert_{M_x \setminus N_x}) - m_1} \; dx. 
\end{split}
\end{equation*}
Statement~(i) follows from this by noting that $\mu_{[2]}(N(\mcx)) = m_1^2 \abs{\mcx} \mck(N)$ (see Equation~\eqref{eq:mmreduction}) and using the various spatial homogeneities that were required. Statement~(ii) is obtained likewise, using additionally that $\EE \Lambda(\mcx) = p \, m_1 \abs{\mcx}$.
\end{proof}
\begin{proof}[Proof of Corollary~\ref{cor:maincontracted}]
   From the definition it is clear that the total variation metric is not affected by 
   changes of scale of the state space, so that
   \begin{equation} \label{eq:dtvcontracted}
      \dtv \bigl(\msl(\xi_{\pi} \kappatinv), Po(\mu\obenpi \kappatinv) \bigr) = \dtv
      \bigl(\msl(\xi_{\pi}), Po(\mu\obenpi) \bigr). 
   \end{equation}
   The definition of $d_1$ and the inequality required for $\dzero$ imply that $\done( \varrho_1 
   \kappatinv, \varrho_2 \kappatinv) \leq \done( \varrho_1, \varrho_2 )$ for all $\varrho_1, \varrho_2 \in 
   \mfn$, whence, by Equation~\eqref{eq:krdtwo},
   \begin{equation} \label{eq:dtwocontracted}
      \dtwo \bigl(\msl(\xi_{\pi} \kappatinv), Po(\mu\obenpi \kappatinv) \bigr) \leq \dtwo
      \bigl(\msl(\xi_{\pi}), Po(\mu\obenpi) \bigr). 
   \end{equation}
   With Equations~\eqref{eq:dtvcontracted} and~\eqref{eq:dtwocontracted} it is seen that all the
   results from Theorem~\ref{thm:main} to Corollary~\ref{cor:mainconst} remain correct if we do the proposed 
   replacements; in particular, the upper bounds stated follow immediately from
   Theorem~\ref{thm:main}.
\end{proof}

\section{Applications} \label{sec:applications} 

We study two applications for a fairly general point process $\xi$ here. The first one concerns the thinning of $\xi$ by covering it with an independent Boolean model. This is up to a few technical adjustments the setting that was used in Section~3.3 of~\cite{schumi2}. We present this application in order to illustrate to what degree the results of the current article improve on the main distance bounds in \cite{schumi2}, and give new insight into the high intensity limit behavior. The second application deals with a Mat{\'e}rn type I thinning of $\xi$. We present it as an example where the rather involved $\bbeta$-term is non-zero and can be reasonably simplified. The bound is compared to a result in \cite{xia05}, where the same thinning was considered for the special case that $\xi$ is a Poisson process.

In this whole section we consider a metric $\tdzero$ on $\RD$ that is generated by a norm,
and use notation of the form $\BB(x,r)$ for closed $\tdzero$-balls in $\RD$ and $\BB^c(x,r)$ for their complements. Write furthermore $\BB_{\mcx}(x,r) := \BB(x,r) \cap \mcx$ for the corresponding balls in $\mcx$ and $\bxc(x,r) := \mcx \setminus \BB_{\mcx}(x,r)$.
We call the subset $\mcx$ of $\RD$ \emph{admissible} if it is compact, of positive volume, and has a boundary $\partial \mcx$ that is of volume zero.

\subsection{Thinning by covering with an independent Boolean model} \label{ssec:ibooleanapplic}

The details for this situation are as follows.
\begin{msetting1}
Suppose that $\mcx \subset \RD$ is admissible and that $\xi$ is a point process on $\mcx$
which has a density $f$ with respect to $P_1 := \poisson(\LebD \vert_{\mcx})$ and finite expectation measure $\mu = \mu_1$ with density $u$. Let $q \in [0,1]$, and take $\Xi$ to be a stationary Boolean model (see \cite{skm87}, Section~3.1) on $\RD$ whose grains are $\tdzero$-balls of random but essentially bounded radius, denoting by $l_1 > 0$ the intensity of the germ process and by $R_i \in \linfty$ the radii of the grains (which are i.i.d.).
This means that $\Xi$ takes the form
\begin{equation}
   \Xi = \bigcup_{i=1}^{\infty} \BB(Y_i, R_i),
\end{equation}
where $Y_i$ are the points of a $\poisson(l_1 \LebD)$-process that is independent of $(R_i)_{i \in \NN}$. 
Assume furthermore that $\xi$ and $\Xi$ are independent, and define a retention field by $\pi(\omega, x) := q \hspace*{1pt} {\rm I}[x \not\in \Xi(\omega)]$ for $\omega \in \Omega$ and $x \in \mcx$. Thinning with respect to $\pi$ corresponds to deleting all the points that are covered by $\Xi$, while retaining uncovered points independently of one another with probability~$q$. $\Diamond$
\end{msetting1}

We aim at applying Corollary~\ref{cor:main} in this setting. Assume without loss of generality that $\PP[R_1 > 0] > 0$ (otherwise Proposition~\ref{prop:iboolean} below is easily checked directly),
and remove from $\Xi(\omega) \cap \mcx$ any lower-dimensional parts, stemming either from balls with radius zero or from balls that only just touch $\mcx$ from the outside, by taking the closure of its interior in $\RD$.
Note that this does not alter the distribution of the obtained thinning, because only a set of volume zero is removed in this way and because $\xi$ and $\pi$ are independent. As a consequence of Proposition~\ref{prop:leps}(iv), where $\mcy = \mcx$ and $\Sigma = \QQ^D \cap \mcx$, we obtain then that $\pi$ has an evaluable path space. Let $N_x := \BB_{\mcx}(x,\bar{r})$ for some $\rbar \geq 2 \linftynorm{R_1}$ and every $x \in \mcx$,
which implies independence of $\pi(x)$ and $\pi \vert_{N_x^c}$ and hence that we can choose $\bgamma \equiv 0$ in Inequality~\eqref{eq:bgammai}.
We set furthermore $r := \ldnorm{R_1}$, so that $r^D = \EE(R_1^D)$.
By the fact that the capacity functional $T_{\Xi}$ of the Boolean model $\Xi$ is given by
\begin{equation*}
   T_{\Xi}(C) := \PP[ \Xi \cap C \neq \emptyset] = 1 - \exp \bigl( - l_1 \, \EE \bigl( \LebD(
   \BB(0,R_1) + C) \bigr) \bigr)
\end{equation*}
for any compact set $C \subset \RD$ (see \cite{skm87}, Equation~(3.1.2)), we obtain for the expectations in the upper bound of Corollary~\ref{cor:main}
\begin{equation}
   \EE \pi(x) = q \bigl( 1 - T_{\Xi}(\{x\}) \bigr) = q e^{- l_1 \EE \abs{\BB(0,R_1)}}
\end{equation}
and
\begin{equation}
   \EE \bigl( \pi(x) \pi(\tx) \bigr) = q^2 \bigl( 1 - T_{\Xi}(\{ x, \tx \}) \bigr) = q^2
   e^{- l_1 \EE \abs{\BB(0,R_1) \cup \BB(\tx-x,R_1)}}.
\end{equation}
As earlier, we use absolute value bars for a measurable subset of $\RD$ to denote its Lebesgue mass. 
Defining $\alpha_D := \abs{\BB(0,1)}$ and $b: \RD \to [0,1]$ by $b(y) := \EE \abs{\BB(0,R_1) \setminus \BB(y,R_1)} \big/ \EE \abs{\BB(0,R_1)}$, we then have the following result.
\begin{prop}  \label{prop:iboolean}
   Under Model Setting 1 laid down above and letting $\rbar \geq 2 \linftynorm{R_1}$, $r := \ldnorm{R_1}$, and
   $N_{\rbar}(\mcx) := \{ (x, \tx) \in \mcx^2; \, \tdzero(x, \tx) \leq \rbar \}$, we obtain that
   \begin{equation*}
   \begin{split}
      \dtv \bigl( \msl(\xi_{\pi}), \poisson(\mu\obenpi) \bigr) &\leq q^2 e^{-2 l_1 \alpha_D r^D} \mu_1^2
         \bigl( N_{\rbar}(\mcx) \bigr)  \\
      &\hspace*{6mm} {} + q^2 \int_{N_{\rbar}(\mcx)} e^{ 
         -l_1 (1+b(\tx-x)) \alpha_D r^D} \; \mu_{[2]} \bigl( d(x,\tx) \bigr)
         \\[0.5mm]
      &\hspace*{6mm} {} + 2 q e^{- l_1 \alpha_D r^D} \abs{\mcx} \bbetasup_{\rbar}, 
   \end{split}
   \end{equation*}
   where $\mu\obenpi = q e^{-l_1 \alpha_D r^D} \mu_1$ and
   \begin{equation} \label{eq:bbetasup}
      \bbetasup_{\rbar} := \sup_{x \in \mcx} \sup_{C \in \mcn(\bxc(x,\rbar))}
         \biggabs{\int_C \bigl[ f_{\bxc(x,\rbar) \cup \{x\}}(\varrho +
         \delta_x) - f_{\bxc(x,\rbar)}(\varrho) u(x) \bigr] \; P_{1,\bxc(x,\rbar)}(d \varrho)}.
   \end{equation}
   \vspace*{-5.5mm}

   \hfill \qed
\end{prop}
\begin{rem}  \label{rem:ibooleancomp1}
   Under Assumption~1 made for Proposition~3.G in \cite{schumi2}, the above estimate can be bounded
   by ${\rm const} \hspace*{-1.2pt} \cdot \hspace*{-1.2pt} \bigl( \rbar^D \abs{\mcx} q^2 e^{- l_1 \alpha_D r^D} 
   + \abs{\mcx} q e^{- l_1 
   \alpha_D r^D} \bbetasup_{\rbar} \bigr)$, which makes it somewhat better than the one in 
   Proposition~3.G (if we accept \smash{$\bbetasup_{\rbar}$} as a natural substitute for $\bbetaind(m)$ in
   \cite{schumi2} and apply Equation~\eqref{eq:dtvcontracted}), also since the result is
   formulated in the stronger $\dtv$-metric instead of~$\dtwo$.

   \nin However, the main point worth noting here is that the derivation above is considerably
   simpler and more elegant than the one for Proposition~3.6, because we do not have to worry about
   covering discretization cuboids. For the same reason we are easily able to work with balls
   that are based on other metrics than the Euclidean one and can write down the explicit constants
   in the upper bound.
\end{rem}
If we assume that $\xi$ is second order stationary (i.e.\ the restriction to $\mcx$ of a second order stationary point process $\zeta$ on $\RD$),
the rather complex second factorial moment measure can be replaced by a term involving the corresponding reduced moment measure. Second order stationarity means that the second moment measure $\mu_2$ of $\zeta$ is locally finite ($\mu_2(B) < \infty$ for every bounded measurable $B \subset \RD$) and invariant under translations along the diagonal $\{(x,x) ; \, x \in \RD \}$ (see~\cite{dvj88}, Definition 10.4.I), and implies stationarity of the expectation measure, so that $\mu_1 = m_1 \LebD$ for some $m_1 \in \Rplus$. It follows from Lemma~10.4.III in~\cite{dvj88} that there is a measure $\mck$ on $\RD$ (unique if $m_1>0$) such that 
\begin{equation}  \label{eq:mmreduction}
   \int_{\RD \times \RD} h(x,\tx) \; \mu_{[2]} \bigl( d(x,\tx) \bigr) = m_1^2 \int_{\RD} \int_{\RD} h(x,x+y)
   \; \mck(dy) \; dx
\end{equation}
for every measurable function $h: \mcx^2 \to \Rplus$. Hence 
\begin{equation}  \label{eq:mmredconcrete}
\begin{split}   
   \int_{N_{\rbar}(\mcx)} e^{-l_1 (1+b(\tx-x)) \alpha_D r^D} \;
      \mu_{[2]} \bigl( d(x,\tx) \bigr) &= m_1^2 \int_{\mcx} \int_{\BB_{\mcx-x}(0,\rbar)} e^{-l_1
      [1+b(y)] \alpha_D r^D} \; \mck(dy) \; dx \\
   &\leq m_1^2 \abs{\mcx} \int_{\BB(0,\rbar)} e^{-l_1 [1+b(y)] \alpha_D r^D} \; \mck(dy). 
\end{split}
\end{equation}
If $\zeta$ is stationary and $m_1>0$,
it can be seen by Equation~\eqref{eq:cm1} (see \cite{skm87}, beginning of Section~4.5) that the measure $\mck$ is given as
\begin{equation}  \label{eq:kmeasure}
   \mck(B) = \frac{1}{m_1} \EE \zeta_0^{\hspace*{0.6pt}!}(B)
\end{equation}
for every Borel set $B \subset \RD$, where $\zeta_0^{\hspace*{0.6pt}!}$ denotes the reduced Palm process of $\zeta$ given a point in~$0$ (see \cite{kallenberg86}, Lemma~10.2 and Section~12.3), so that
$\EE \zeta_0^{\hspace*{0.6pt}!}(B)$ can be interpreted as the expected number of points of $\zeta$ in $B$ given there is a point in $0$. The measure $\mck$ is usually referred to as \emph{second reduced moment measure}, although some authors prefer  defining it as $m_1$ (or even $m_1^2$) times the above measure.
Set furthermore $K(\tr) := \mck \bigl( \BB(0,\tr) \bigr)$ for every $\tr \in \Rplus$, which, if $\tdzero$ is the Euclidean metric, defines Ripley's $K$-function. 

We examine the situation of Corollary~3.H in \cite{schumi2}, waiving two technical conditions that were needed there, but insisting on second-order stationarity in order to bring the second summand in the upper bound in a nicer form.

\begin{msetting1p}
Suppose that $J \subset \RD$ is admissible, that $n \in \NN$, and that $\mcx = \mcx_n = \kappa_n^{-1}(J)$, where $\kappa_n(x) = (1/n)x$ for every $x \in \RD$. Let $\xi$ be a second order stationary point process on~$\mcx$  which has density $f$ with respect to $P_1 = \poisson(\LebD \vert_{\mcx})$ and expectation measure $\mu = \mu_1 = m_1 \LebD \vert_{\mcx}$ for some $m_1 \in \Rplus$. We assume that $\xi$ is the restriction to $\mcx$ of one and the same point process $\zeta$ on $\RD$ for every $n$, and suppress the index $n$ in any quantities that depend on $n$ only by virtue of this restriction.
Choose a sequence $(q_n)_{n \in \NN}$ with $1/n^D \leq q_n \leq 1$, and a sequence $(\Xi_n)_{n \in \NN}$ of stationary Boolean models on $\RD$ of $\tdzero$-balls with radii $R_{n,i} \in \linfty$ (i.i.d.\ for every $n$) and germ process intensity $l_1 > 0$ such that
\begin{equation}
   r_n := \bigldnorm{R_{n,1}} = \Bigl( \frac{1}{l_1 \alpha_D} \log(q_n n^D) \Bigr)^{1/D}.
\end{equation}
Assume that $\xi$ and $\Xi_n$ are independent for every $n \in \NN$, and define retention fields by $\pi_n(\omega, x) := q_n \hspace*{0.5pt} {\rm I}[x \not\in \Xi_n(\omega)]$ for $\omega \in \Omega$ and $x \in \mcx$. Let furthermore $\rbar_n \geq 2 \linftynorm{R_{n,1}}$, and note that $\mu\obenpi \kappa_n^{-1} = \frac{1}{n^D} \mu_1 \kappa_n^{-1} = m_1 \LebD \vert_J$. $\Diamond$
\end{msetting1p}

By Inequality~\eqref{eq:mmredconcrete} and as $\dtv \bigl( \msl(\xi_{\pi_n} \kappa_n^{-1}), \poisson(m_1 \LebD \vert_J) \bigr) = \dtv \bigl( \msl(\xi_{\pi_n}), \poisson(\frac{m_1}{n^D} \LebD \vert_{\mcx}) \bigr)$ by
Equation~\eqref{eq:dtvcontracted}, we have the following consequence of Proposition~\ref{prop:iboolean}.
\begin{cor}  \label{cor:iboolean}
   Under Model Setting 1$'$ we obtain that  
   \begin{equation*}
   \begin{split}
      \dtv \bigl( \msl(&\xi_{\pi_n} \kappa_n^{-1}), \poisson(m_1 \LebD \vert_J) \bigr) \\[1mm]
      &\leq m_1^2 \abs{J} \alpha_D \Bigl( \frac{\rbar_n}{n} \Bigr)^D + m_1^2 \abs{J} q_n \int_{\BB(0,
         \rbar_n)} (q_n n^D)^{-b(y)} \; \mck(dy) + 2 \abs{J}
         \bbetasup_{\rbar_n} \\
      &\leq \abs{J} \biggl( m_1^2 \alpha_D \Bigl( \frac{\rbar_n}{n} \Bigr)^D + m_1^2 q_n K(\rbar_n) +
         2 \bbetasup_{\rbar_n} \biggr),
   \end{split}
   \end{equation*}
   where $\bbetasup_{\rbar}$ was defined in Equation~\eqref{eq:bbetasup}. \hfill \qed
\end{cor}
\begin{rem}
   Note that Assumption 1b) in \cite{schumi2} implies that $K(\rbar_n) = O(\rbar_n^D)$ for $n \to \infty$.
   Compared with the corresponding result in \cite{schumi2} (Corollary~3.H) we therefore have again somewhat 
   better bounds
   with explicit constants that were obtained in a more direct way.
\end{rem}

A rather nice result can be derived from Corollary~\ref{cor:iboolean} in the Poisson case.
\begin{cor}  \label{cor:ibooleanppp}
   Under Model Setting 1$'$ and the additional assumptions that 
   $\zeta$ is a Poisson process with expectation measure $m_1 \LebD$, that $q_n=1$ for every $n \in \NN$, that
   $\tdzero$ is the Euclidean metric, and that $\linftynorm{R_{n,1}} = O(r_n)$ for $n \to \infty$,
   we have
   \begin{equation*}
      \dtv \bigl( \msl(\xi_{\pi_n} \kappa_n^{-1}), \poisson(m_1 \LebD \vert_J) \bigr) 
      = O \bigl( (\log n)^{-(D-1)} \bigr) \textfor n \to \infty.
   \end{equation*}
\end{cor}
\begin{rem} \label{rem:ibooleanconv}
   The following (partly heuristical) arguments suggest that the order claimed in
   Corollary~\ref{cor:ibooleanppp} is sharp for $m_1>0$. Assume for simplicity that $J = [0,1]^D$.
   
   \nin
   For $D=1$ it is readily understood that $\xi_{\pi_n} \kappa_n^{-1}$ cannot converge in 
   distribution to a Poisson process. The reason is that the uncovered part of $\RR$ in the domain
   of the contraction, i.e.\ $\RR \setminus \Xi_n$, is made up of intervals whose lengths are exponentially 
   distributed with mean $1/l_1$ (no matter how the $R_{n,i}$ are distributed),
   so that the 
   probability of having two or more points within the first uncovered interval that lies completely in
   $\Rplus$ does not depend on $n$. Hence we have a constant positive probability 
   that the first two points of the thinned contracted point process in $\Rplus$ are within distance $1/
   (l_1 n)$, say, 
   which cannot be true for a sequence that converges towards a homogeneous Poisson process.
   
   \nin
   For $D \geq 2$ the situation is more complicated. By Theorem~1 in \cite{hall85} (compare also
   statement~(ii) on page~244), it can be seen that the uncovered ``chinks'' of $\Xi_n$ in the domain of 
   the contraction have a volume that is of order $1 \big/ \bigl(\log(n^D)\bigr)^{D-1}$ for large $n$,
   so that the argument of the constant-sized chinks is not valid for general $D$.
   Heuristically, the order of the chink volumes together with the fact that the number of chinks in a bounded
   measurable set is Poisson distributed (see \cite{hall85}, p.~244, statement~(i); note the slightly different 
   scaling) suggest that we can think of the process $\xi_{\pi_n} \kappa_n^{-1}$ for $n$ large as a compound 
   Poisson process \raisebox{0pt}[11pt][0pt]{$\sum_{i=1}^{V_n} Z_i\obenn \delta_{S_i\obenn}$} with intensity of 
   the Poisson 
   process \raisebox{0pt}[3pt][5pt]{$\sum_{i=1}^{V_n} \delta_{S_i\obenn}$} of order $\bigl( \log(n^D)
   \bigr)^{D-1}$
   and i.i.d.\ clump sizes $Z_i\obenn$, for which $\PP[Z_1\obenn \geq 1]$ is of order $\bigl( \log(n^D)
   \bigr)^{-(D-1)}$ and $\PP[Z_1\obenn \geq 2]$ is of
   order \raisebox{0pt}[10.5pt][0pt]{$\bigl( \log(n^D) \bigr)^{-2(D-1)}$} (by the fact that $\xi$ is a Poisson
   process).
   It is
   easily seen that such a process converges towards a homogeneous Poisson process $\eta$ as $n \to \infty$ by
   noting that \raisebox{0pt}[12pt][0pt]{$\dtv \bigl(\msl(\sum_{i=1}^{V_n} Z_i\obenn \delta_{S_i\obenn}), \msl
   (\sum_{i=1}^{V_n}
   {\rm I}[Z_i\obenn \geq 1] \delta_{S_i\obenn}) \bigr) \to 0$}, and
   \raisebox{0pt}[12pt][3.5pt]{$\sum_{i=1}^{V_n}
   {\rm I}[Z_i\obenn \geq 1] \delta_{S_i\obenn} \inlawto \eta$}, 
   but that its convergence rate in the total variation metric is bounded from below by $\PP[\exists i \in \{1, 
   \ldots, V_n\} : Z_i\obenn \geq 2]$, which is of order $\bigl( \log(n^D) \bigr)^{-(D-1)}$ or, what is the 
   same, order $(\log n)^{-(D-1)}$.    
\end{rem}
\begin{proof}[Proof of Corollary~\ref{cor:ibooleanppp}]
   Our starting point is the first upper bound in Corollary~\ref{cor:iboolean}, where we set $\rbar_n := 2
   \linftynorm{R_{n,1}}$. Since $\zeta$ is Poisson, the 
   third summand is zero, whereas the first summand is clearly $O \bigl( (\log n)^{-(D-1)} \bigr)$. We 
   investigate the integral in the second summand. We have $\mck = \LebD$ by
   Equation~\eqref{eq:kmeasure} in combination with $\msl(\zeta_0^{\hspace*{0.6pt}!}) = \msl(\zeta)$ (see
   \cite{moellerwaage04}, Proposition~C.2). Define $\tb: [0,2] \to [0,1]$ by $\tb(u) = 
   \abs{\BB(0,1) \setminus \BB(y,1)} \big/ \abs{\BB(0,1)}$, where $y$ is an arbitrary element of $\RD$ with 
   $\abs{y} = u$, and note that \raisebox{0pt}[10pt][0pt]{$b(y) \geq \tb(2\abs{y}/\rbar_n)$}
   for $y \in \BB(0,\rbar_n)$. Since $\tdzero$ is the Euclidean metric, it can be shown that there is a constant
   $\kappa >0$ such that $\tb(u) \geq \frac{\kappa}{2} u$ for every $u \in [0, 2]$. 
   Writing $\omega_D$ for the surface area of the unit sphere in $\RD$, we then can bound the required
   integral as 
   \begin{equation*}
   \begin{split}
      \int_{\BB(0, \rbar_n)} (q_n n^D)^{-b(y)} \; \mck(dy) &\leq \int_{\BB(0, \rbar_n)}
         (n^D)^{-\kappa \abs{y}/\rbar_n} \; \LebD(dy) \\
      &= \frac{\omega_D}{\kappa^D} \rbar_n^D \int_0^{\kappa} (n^D)^{-r} r^{D-1} \; dr \\[0.5mm]
      &\leq \frac{(D-1)! \, \omega_D}{\kappa^D} \rbar_n^D \bigl(\log(n^D)\bigr)^{-D} (1 - n^{- \kappa D}), 
   \end{split}
   \end{equation*}
   where the last inequality follows from multiple integration by parts.
   Since $\rbar_n = O \bigl( (\log(n^D))^{1/D} \bigr)$, we thus obtain that also the second summand in 
   the first upper bound of Corollary~\ref{cor:iboolean} is $O \bigl( (\log n)^{-(D-1)} \bigr)$.
\end{proof}

\subsection{Thinning by Mat{\'e}rn type I competition} \label{ssec:maternapplic}

Again we base our retention field on a random closed set $\Xi$, but this time we choose a situation where $\Xi$ is completely determined by the point process $\xi$. The resulting procedure is the one used for the construction of the Mat{\'e}rn type I hard core process, in which a point is deleted whenever there is any other point within a fixed distance $r$. The details are as follows. 
\begin{msetting2}
Suppose that $r >0$ and that $\mcx, \mcx' \subset \RD$ are two compact sets, where $\mcx$ is admissible and $\BB(\mcx, r) \subset \mcx'$. Let furthermore $\xi$ be a stationary point process on $\mcx'$ (i.e.\ the restriction of a stationary point process $\zeta$ on $\RD$) which has density $f$ with respect to $P_1 := \poisson(\LebD \vert_{\mcx'})$ and a finite expectation measure $\mu = \mu_1 = m_1 \LebD \vert_{\mcx'}$ for some $m_1 \in \Rplus$. By defining $\xi$ on $\mcx'$, but considering the thinned point process only on $\mcx$, we avoid boundary effects, which would lead to more complicated notation because of spatial
inhomogeneities in the thinned process.

\nin
In order to have a $\Xi$ whose realizations are closed sets that are jointly separable, we proceed as follows. Write $\xi$ as $\sum_{i=1}^V \delta_{S_i}$, where $V$ and $S_i$ are $\sigma(\xi)$-measurable random elements with values in $\Zplus$ and $\mcx'$, respectively, and denote by $T_i$ the $\tdzero$-distance between $S_i$ and its nearest neighbor. Let then
\begin{equation*}
   B_i(\omega) := \bigl\{ y \in \mcx'; \, \tfrac13 \min(T_i(\omega),r) \leq \tdzero(y,S_i(\omega)) \leq r
   \bigr\}
\end{equation*}
for every $\omega \in \Omega$, and set
\begin{equation*}
   \Xi := \bigcup_{i=1}^{\infty} B_i.
\end{equation*}
Choose $q \in [0,1]$ and
define a retention field on $\mcx'$ by setting $\pi(\omega,x) := q \hspace*{1pt} {\rm I}[x \not\in \Xi(\omega)]$ if $\omega \in \Omega$ and $x \in \mcx$ and $\pi(\omega,x) := 0$ if $\omega \in \Omega$ and $x \in \mcx' \setminus \mcx$.
Note that
\begin{equation} \label{eq:maternpi}
   \pi(\omega,s) = q {\rm I} \bigl[ s \in \mcx, \xi(\omega)(\dot{\BB}(s,r))=0 \bigr]
\end{equation}
for $\omega \in \Omega$ and $s \in \xi(\omega)$, where $\dot{\BB}(x,r) := \BB(x,r) \setminus \{x\}$ for every $x \in \mcx$. Hence, on $\mcx$, thinning with respect to $\pi$ corresponds to deleting all those points that see any other point of the
process within distance $r$ (regardless whether this point is itself deleted or not), while retaining points that do not have this property independently of one another with probability $q$. $\Diamond$
\end{msetting2}

This time, we aim at applying Theorem~\ref{thm:main} for the state space $\mcx'$. By Proposition~\ref{prop:leps}(iv), $\pi$ has an evaluable path space (after removing from $\Xi(\omega) \cap \mcx$ possible lower-dimensional parts by taking the closure of its interior in $\RD$, which has no influence on the distribution of the resulting thinning).
Since $\pi$ is completely determined by $\xi$ we have the corresponding Dirac measure as a regular conditional distribution of $\pi$ given the value of $\xi$. Condition~\eqref{eq:locdepcond} is satisfied for a catchment radius of $R=r$, so that $N_x := \BB_{\mcx'}(x, \rbar)$ for some $\rbar \geq 2r$ is a legitimate choice for the neighborhoods of strong dependence. We can furthermore choose $\bgamma \equiv 0$ in Inequality~\eqref{eq:bgamma}.

Write $\xi_x^{\hspace*{0.6pt}!}$ for the reduced Palm process of $\xi$ given a point in $x$, and $\xi_{x,\tx}^{\hspace*{0.6pt}!}$ for the second-order reduced Palm process of $\xi$ given points in $x$ and $\tx$ (see \cite{kallenberg86}, Section~12.3, pp.\ 109 \& 110; note that $\nu_n' = \mu_{[n]}$ for obtaining the distributions of the $n$-th order reduced Palm processes). The first and second order Campbell-Mecke equations state that
\begin{equation} \label{eq:cm1}
   \EE \biggl( \int_{\mcx'} h(x, \xi - \delta_x) \; \xi(dx) \biggr) = \int_{\mcx'} \EE h(x, 
   \xi_x^{\hspace*{0.6pt}!}) \; \mu_1(dx)  
\end{equation}
and
\begin{equation} \label{eq:cm2}
   \EE \biggl( \int_{\mcxpsq} h(x, \tx, \xi - \delta_x -\delta_{\tx}) \; \xi^{[2]} \bigl( d(x,\tx) \bigr) 
   \biggr) = \int_{\mcxpsq} \EE h(x, \tx, \xi_{x,\tx}^{\hspace*{0.6pt}!}) \; \mu_{[2]} \bigl( d(x,\tx) \bigr)  
\end{equation}
for non-negative measurable functions $h$.
These equations follow immediately by standard extension arguments from the defining equations of Palm
processes (see e.g.\ \cite{moellerwaage04}, Equation~(C.4), for the first one).
We then obtain by Equation~\eqref{eq:maternpi} that
\begin{equation} \label{eq:maternLambda1}
\begin{split}
   \EE \Lambda(A) &= \EE \biggl( \int_{A \cap \mcx} q {\rm I} \bigl[ \xi(\BBdot(x,r)) = 0 \bigr]
      \; \xi (dx) \biggr) \\
   &= q \int_{A \cap \mcx} \PP \bigl[ \xi_x^{\hspace*{0.6pt}!} \bigl( \BB(x,r) \bigr) = 0 \bigr] \; \mu_1(dx)
      = m_1 q \bigl( 1-G(r) \bigr) \abs{A \cap \mcx}
\end{split}
\end{equation}
for any $A \in \mcb = \mcb(\mcx')$, where $G: \Rplus \to [0,1]$, $G(\tr) := \PP \bigl[ \zeta_{0}^{\hspace*{0.6pt}!} \bigl( \BB(0,\tr) \bigr) \geq 1 \bigr] = \PP \bigl[ \zeta_{x}^{\hspace*{0.6pt}!} \bigl( \BB(x,\tr) \bigr) \geq 1 \bigr]$ for arbitrary $x \in \RD$, denotes the nearest neighbor function of $\zeta$, i.e.\ the distribution function of the distance from a ``typical point'' to its nearest neighbor, which is a frequently used tool in spatial statistics; see e.g.~\cite{baddeley07} (Section~3.4), \cite{diggle03}, or \cite{moellerwaage04}. In a very similar way, using in addition Equation~\eqref{eq:mmreduction} to obtain the last equality, we have with $N_{\rbar}(\mcx') := \{ (x, \tx) \in \mcxpsq; \, \tdzero(x, \tx) \leq \rbar \}$ and $N_{\rbar}(\mcx) := N_{\rbar}(\mcx') \cap \mcx^2$ that 
\begin{equation} \label{maternLambda2}
\begin{split}
   \EE \Lambda^{[2]} \bigl( N_{\rbar}(\mcx') \bigr) &= \EE \biggl( \int_{N_{\rbar}(\mcx)} q^2 {\rm I} 
      \bigl[ \xi(\BBdot(x,r)) = 0 \bigr] {\rm I} \bigl[ \xi(\BBdot(\tx,r)) = 0 \bigr] \; \xi^{[2]}
      \bigl( d(x,\tx) \bigr) \biggr) \\
   &= q^2 \int_{N_{\rbar}(\mcx)} {\rm I}\bigl[ \tdzero(x,\tx) > r \bigr] \PP \bigl[ \xi_{x,
      \tx}^{\hspace*{0.6pt}!} \bigl(
      \BB(x,r) \cup \BB(\tx,r) \bigr) = 0 \bigr] \; \mu_{[2]} \bigl( d(x,\tx) \bigr) \\
   &= m_1^2 q^2 \int_{\mcx} \int_{(\BB(0,\rbar) \setminus \BB(0,r)) \cap (\mcx-x)} \bigl( 1-G_{2,y}(r) 
      \bigr) \; \mck(dy) \; dx,
\end{split}
\end{equation}
where $G_{2,y}: \Rplus \to [0,1]$, $G_{2,y}(\tr) := \PP \bigl[ \zeta_{0,y}^{\hspace*{0.6pt}!} \bigl( \BB(0,\tr) \cup \BB(y,\tr) \bigr) \geq 1 \bigr] = \PP \bigl[ \zeta_{x,x+y}^{\hspace*{0.6pt}!} \bigl( \BB(x,\tr) \cup \BB(x+y,\tr) \bigr) \geq 1 \bigr]$ for arbitrary $x \in \RD$, are the natural two-point analogs of the $G$-function, with $y \in \RD$.

Finally, the term $\bbeta(x) = \bbeta_{\rbar}(x)$ is zero for $x \in \mcx' \setminus \mcx$, and can be evaluated for $x \in \mcx$ as
\begin{equation} \label{eq:maternbeta1}
\begin{split}
   \bbeta_{\rbar}(x) &= \int_{\mfn(\BB(x,r))} q {\rm I}[\varrhoint = 0] \int_{\mfn(\bxpc(x,
      \rbar-r))} \bigabs{\fbar(\varrhoext, \varrhoint + \delta_x)} \; P_{1, \bxpc(x,\rbar-r)}(d
      \varrhoext) \; P_{1, \BB(x,r)}(d\varrhoint) \\
   &= q e^{-\alpha_D r^D} \int_{\mfn(\bxpc(x, \rbar-r))} \bigabs{\fbar(\varrhoext, \delta_x)} \;
      P_{1, \bxpc(x,\rbar-r)}(d \varrhoext) \\[1mm]
   &= q e^{-\alpha_D r^D} \EE \Bigabs{f_{\BB(x,r) \mvert \bxpc(x,\rbar-r)}(\delta_x; \xi \vert_{\bxpc(x,
      \rbar-r)}) - f_{\BB(x,r)}(\delta_x)},
\end{split}
\end{equation}
where we set $f_{\Aint \mvert \Aext}(\varrhoint; \varrhoext) := f_{\Aext \cup \Aint}(\varrhoext+\varrhoint) \big/ f_{\Aext}(\varrhoext)$ if $f_{\Aext}(\varrhoext) > 0$ and $f_{\Aint \mvert \Aext}(\varrhoint; \varrhoext) := 0$ otherwise.

The expectation in the last line of \eqref{eq:maternbeta1} is much simpler than the general expression we have for $\bbeta$, but is typically still hard to estimate. A more directly applicable estimate, which, however, is very rough and works only with point processes that are not too extreme in a certain sense, is given as follows. Choose $\rbar := 2r$, and assume that $f > 0$ $P_1$-almost surely.
We then obtain with $\etaext \sim P_{1, \bxpc(x,r)}$ that
\begin{equation} \label{eq:maternbeta2}
\begin{split}
   \bbeta(x) &= \int_{\mfn(\BB(x,r))} q {\rm I}[\varrhoint = 0] \int_{\mfn(\bxpc(x,r))}
      \bigabs{\fbar(\varrhoext, \varrhoint + \delta_x)} \; P_{1, \bxpc(x,r)}(d\varrhoext) \;
      P_{1, \BB(x,r)}(d\varrhoint) \hspace*{-30mm} \\
   &= q \int_{\mfn} {\rm I}[\sigma \vert_{\BB(x,r)} = 0] \Bigabs{ f(\sigma \vert_{\bxpc(x,r)} + 
      \delta_x) - f_{\bxpc(x,r)}(\sigma \vert_{\bxpc(x,r)}) f_{\BB(x,r)}(\delta_x) } \; 
      P_1(d\sigma) \\
   &= q \int_{\mfn} {\rm I}[\sigma(\BB(x,r)) = 0] \biggabs{ 1 - \frac{f_{\bxpc(x,r)}(\sigma
      \vert_{\bxpc(x,r)}) f_{\BB(x,r)}(\delta_x)}{f(\sigma \vert_{\bxpc(x,r)} +
      \delta_x)} } f(\sigma + \delta_x) \; P_1(d \sigma) \\
   &\leq m_1 q \bigl( 1 - G(r) \bigr) \bigglinftynorm{1 - \frac{f_{\bxpc(x,r)}(\etaext)
      f_{\BB(x,r)} (\delta_x)}{f(\etaext + \delta_x)}}
\end{split}
\end{equation}
for $x \in \mcx$,
where we used in the last line that
\begin{equation*}
   m_1 \hspace*{0.5pt} \EE \bigl( h(x, \xi_x^{\hspace*{0.6pt}!}) \bigr) = \int_{\mfn} h(x, \sigma) f(\sigma + \delta_x) \;
   P_1(d \sigma)
\end{equation*}
for almost every $x \in \mcx'$ and every non-negative measurable function $h$, which is a consequence of the first Campbell-Mecke equation and of Equation~\eqref{eq:pregnz} with $N_x = \{x\}$.
We assume that the $\linftynorm{\cdot}$-term is bounded by a constant $M \in \Rplus$ (that depends neither on $r$ nor on $x$).
Two examples where this is satisfied are given at the end of this subsection.

Pluging \eqref{eq:maternLambda1} to \eqref{eq:maternbeta2} into Theorem~\ref{thm:main}(i) and choosing now $\rbar := 2r$ everywhere yields the following result. 
\begin{prop}  \label{prop:matern}
   Under Model Setting 2 laid down above, we obtain that
   \begin{equation*}
   \begin{split}
      \dtv \bigl( \msl( \xi_{\pi}), \poisson(\mu\obenpi) \bigr)
      &\leq m_1^2 \abs{\mcx} 2^D \alpha_D r^D q^2 \bigl( 1 - 
         G(r) \bigr)^2 + m_1^2 \abs{\mcx} q^2 \hspace*{-4pt} \int_{\BB(0,2r) \setminus \BB(0,r)}
         \hspace*{-4pt}\bigl( 1-G_{2,y}(r) \bigr) \; \mck(dy) \\[-0.5mm]
      &\hspace*{6mm} {} + q e^{-\alpha_D r^D} \int_{\mcx} \EE \Bigabs{f_{\BB(x,r) \mvert
         \bxpc(x,r)} (\delta_x; \xi \vert_{\bxpc(x,r)}) - f_{\BB(x,r)}(\delta_x)} \; dx,
   \end{split}
   \end{equation*}
   where $\mu\obenpi = m_1 q (1-G(r)) \LebD \vert_{\mcx}$.
   If in addition $f > 0$ $P_1$-almost surely and the $\linftynorm{\cdot}$-term in
   Inequality~\eqref{eq:maternbeta2} is uniformly bounded by $M \in \Rplus$, then the last summand can be 
   estimated further by
   \begin{equation*}
      m_1 \abs{\mcx} q \bigl( 1 - G(r) \bigr) M.
   \end{equation*}
   \vspace*{-8.5mm}

   \hfill \qed
\end{prop}
\begin{rem}
   In order to obtain an integrand that does not depend on $x$ in the last summand of the above bound, we
   can either proceed as in Corollary~\ref{cor:mainconst}, applying Inequality~\eqref{eq:extendingg}
   and as a consequence replace $\bxpcstd(x,r)$ by $M_x \setminus \BB(x,r)$ for bounded ``outer neighborhoods'' 
   $M_x$ that are shifted copies of one another and that all contain the set $\mcx'$; or we can proceed as in 
   Remark~\ref{rem:mainconst}, using a Gibbs construction on the whole of $\RD$ and as a consequence replace 
   $\bxpcstd(x,r)$ by $\RD \setminus \BB(x,r)$.   
\end{rem}

If $\xi$ is a homogeneous Poisson process and $q=1$, then $\xi_{\pi}$ is the usual Mat{\'e}rn type I hard core process restricted to $\mcx$, and the above bound takes especially simple form.
\begin{cor} \label{cor:matern}
   Under Model Setting 2 and the additional assumptions that $\zeta$ is a Poisson process and $q=1$, we
   have
   \begin{equation*}
   \begin{split}
      \dtv \bigl( \msl( \xi_{\pi}), \poisson(l_1 \LebD \vert_{\mcx} ) \bigr)
         &\leq \abs{\mcx} 2^D \alpha_D r^D l_1^2 + m_1^2 \abs{\mcx} \int_{\BB(0,2r)
         \setminus \BB(0,r)} e^{-m_1 \abs{\BB(0,r) \cup \BB(y,r)}} \; dy \\
      &\leq \abs{\mcx} 2^D \alpha_D r^D l_1^2 (1 + e^{\frac{1}{2} m_1 \alpha_D r^D}),
   \end{split}
   \end{equation*}
   where $l_1 := m_1 e^{-m_1 \alpha_D r^D}$.
\end{cor}
\begin{proof}
   The first inequality follows directly from Propositon~\ref{prop:matern} by the fact that $\zeta$ is a Poisson 
   process, and hence the last summand is zero and the one- and two-point $G$-functions can be easily
   computed by $\msl(\zeta_0^{\hspace*{0.6pt}!}) = \msl(\zeta)$ and $\msl(\zeta_{0,y}^{\hspace*{0.6pt}!}) =
   \msl(\zeta)$ (see \cite{moellerwaage04}, Proposition~C.2, for the one-point case; the two-point case is
   proved in the analogous way).
   
   The second inequality is a consequence of $\abs{\BB(0,r) \cup \BB(y,r)} \geq \frac32 \abs{\BB(0,r)}$ for
   $\abs{y} \geq r$. The latter is due to the fact that all the $\tdzero$-balls of fixed radius are
   translated copies of one another which are convex and symmetric with respect to their centers, and can be
   seen as follows.
   The symmetry implies that any hyperplane through the origin divides
   the volume of $\BB(0,r)$ in half, while the convexity implies the existence of a supporting hyperplane $H_x$ 
   at every point $x$ of the boundary of $\BB(0,r)$, which means that $H_x$ contains $x$ and that $\BB(0,r)$ 
   lies completely in one of the closed half-spaces defined by $H_x$. Thus $H_x - x$ and $H_x$ divide $\RD$
   into three parts, each of which contains half of $\BB(0,r)$ or $\BB(x,r)$, whence we obtain
   that $\abs{\BB(0,r) \cup \BB(x,r)} \geq \frac32 \abs{\BB(0,r)}$ for $\abs{x} = r$. Clearly,
   $\abs{\BB(0,r) \cup \BB(y,r)} \geq \abs{\BB(0,r) \cup \BB(x,r)}$ if $\abs{y} \geq \abs{x}$.
\end{proof}
\begin{rem}
   In Theorem~6.6 of \cite{xia05} a situation very similar to the one in
   Corollary~\ref{cor:matern} for the special case that we choose $\tdzero$ to be the Euclidean metric was 
   considered. 
   The only substantial difference is
   that in \cite{xia05} the Poisson process $\xi$ is defined on $\mcx$ instead of the superset $\mcx'$,
   which leads to less competition near the boundary of $\mcx$ and consequently to a non-stationary thinned 
   process. However, this difference enters the upper bounds in \cite{xia05} only insofar as balls
   are always restricted to $\mcx$ instead of being balls in $\RD$.
   
   \nin
   Disregarding these boundary effects, we see that the estimates in Corollary~\ref{cor:matern} are slightly 
   better than the one for the total variation in
   \cite{xia05}, because our second estimate above is bounded by $2 \abs{\mcx} m_1 \alpha_D
   (2r)^D l_1$, which is exactly the estimate in \cite{xia05} if we adapt it to our notation.
   
   \nin
   The main reason for formulating Corollary~\ref{cor:matern}, however, was not to improve on this earlier 
   bound, but to demonstrate that Proposition~\ref{prop:matern}, which holds for a much greater class of point 
   processes, provides a reasonable estimate in the special case of a Poisson process.
\end{rem}

We end this subsection by giving two examples of point processes for which the additional boundedness condition in Proposition~\ref{prop:matern} is satisfied. 
\begin{example}
   Consider a point process density of the form $f(\varrho) = \lambda^{\abs{\varrho}}    
   \tf(\varrho)$ for a function $\tf$ that is bounded and bounded away from zero.
   One particular instance of such a density is given as
   \begin{equation*}
      f(\varrho) := \kappa \lambda^{\abs{\varrho}} \Bigl( {\rm I} \Bigl[ \min_{\hspace*{0.3pt} s, \ts \in 
      \varrho, \hspace*{0.8pt} s \neq \ts} \tdzero(s,\ts) \leq \tr \Bigr]
      + \gamma {\rm I} \Bigl[ \min_{\hspace*{0.3pt} s, \ts \in \varrho, \hspace*{0.8pt} s \neq \ts}
      \tdzero(s,\ts) > \tr \Bigr] \Bigr)
   \end{equation*}
   for every $\varrho \in \mfn = \mfn(\mcx')$, where $\gamma, \tr, \lambda > 0$ are parameters and $\kappa > 0$ 
   is a normalizing constant. Note that for $\gamma = 1$ we obtain the density of the
   $\poisson(\lambda \LebD \vert_{\mcx'})$-process.
\end{example}
\begin{example}[Strauss process]
   Consider the Strauss process with range
   $\tr \in [0,r]$ and interaction parameter $\gamma \in (0,1]$ (see \cite{moellerwaage04}, Section~6.2.2).
   This process has a density given by
   \begin{equation*}
      f(\varrho) := \kappa \lambda^{\abs{\varrho}} \gamma^{c_{\tr}(\varrho)}
   \end{equation*}
   for every $\varrho \in \mfn$,
   where $c_{\tr}(\varrho) := \sum_{s,\ts \in \varrho, s \neq \ts} {\rm I}[\tdzero(s,\ts) \leq \tr]$ counts the
   pairs of points that lie within distance $\tr$ of one another, $\lambda > 0$ is an intensity
   parameter, and $\kappa >0$ is a normalizing constant. Then, for~$x \in \mcx$ and
   $\sigmaext \in \mfn \bigl( \bxpcstd(x,r) \bigr)$,
   \begin{equation*}
   \begin{split}
      &\frac{f_{\bxpc(x,r)}(\sigmaext) f_{\BB(x,r)}(\delta_x)}{f(\sigmaext +
         \delta_x)} \\[0.5mm]
      &\hspace*{3mm} \leq  \frac{1}{\kappa \hspace*{1pt} \lambda^{\abs{\sigmaext} + 1 
         } \hspace*{1pt} \gamma^{c_{\tr}(\sigmaext)}} \int_{\mfn(\BB(x,r))} \kappa   
         \lambda^{\abs{\sigmaext} + \abs{\varrhoint}} \gamma^{c_{\tr}(\sigmaext)} 
         \; P_{1,\BB(x,r)}(d\varrhoint) \\[-1.5mm]
      &\hspace*{62mm} \cdot \int_{\mfn(\bxpc(x,r))} \kappa \lambda^{\abs{\varrhoext}+1}
         \; P_{1,\bxpc(x,r)}(d\varrhoext) \\[0.5mm]
      &\hspace*{3mm} = \kappa \int_{\mfn(\BB(x,r))} \int_{\mfn(\bxpc(x,r))} \lambda^{\abs{\varrhoint}} \,
         \lambda^{\abs{\varrhoext}} \;
         P_{1,\bxpc(x,r)}(d\varrhoext) \; P_{1,\BB(x,r)}(d\varrhoint) \\[0.5mm]
      &\hspace*{3mm} = \kappa \int_{\mfn} \lambda^{\abs{\sigma}} \; P_1(d\sigma) = \kappa e^{\lambda-1},
   \end{split}
   \end{equation*}
   where we used for the inequality that $c_{\tr}(\sigma + \varrho) \geq c_{\tr}(\sigma)$ for all $\varrho, 
   \sigma \in \mfn$ and $c_{\tr}(\sigma + \delta_x) = c_{\tr}(\sigma)$ for $\sigma \in \mfn$ and $x \in \mcx$ 
   with $\sigma(\BB_{\mcx'}(x, \tr)) = 0$.
   Thus, $M$ may be chosen to be $\max(1, \kappa e^{\lambda-1}-1)$.
\end{example}

\section*{Appendix}

\renewcommand{\thesubsection}{\Alph{section}.\arabic{subsection}}
\renewcommand{\theprop}{\Alph{section}.\Alph{prop}}
\renewcommand{\theequation}{\Alph{section}.\arabic{equation}}
\setcounter{section}{1}
\setcounter{subsection}{0}
\setcounter{prop}{0}
\setcounter{equation}{0}

In what follows we formulate and prove some of the more technical results needed in the main part of this article.

\subsection{Density formulae used for Theorem~\ref{thm:bb}}

\begin{prop} \label{prop:pregnz}
For a point process $\xi$ on $\mcx$ with density $f$ with respect to $P_1$ and finite expectation measure, and for a neighborhood structure $(N_x)_{x \in \mcx}$ that satisfies Condition~\eqref{eq:nhcondition}, we have
\begin{equation*}
   \EE \biggl( \int_{\mcx} h( x, \xi \vert_{N_x^c}) \; \xi(dx) \biggr) = \int_{\mcx} \int_{\mfn(N_x^c)} h
   (x,\varrho) f_{N_x^c \cup \{x\}}(\varrho + \delta_x) \; P_{1,N_x^c}(d \varrho) \; \alpha(dx) 
\end{equation*}
for every non-negative or bounded measurable function $h : \mcx \times \mfn \to \RR$.
\end{prop}
\begin{cor}[Generalized Nguyen-Zessin formula on compact spaces] \label{cor:gnz}
   Suppose that the conditions of Proposition~\ref{prop:pregnz} hold and that $f$ is 
   hereditary. We then have
   \begin{equation*}
      \EE \biggl( \int_{\mcx} h( x, \xi \vert_{N_x^c}) \; \xi(dx) \biggr) = \int_{\mcx} \EE \bigl(
      h( x, \xi \vert_{N_x^c}) g(x, \xi \vert_{N_x^c}) \bigr) \; \alpha(dx)
   \end{equation*}
   for every non-negative or bounded measurable function $h : \mcx \times \mfn \to \RR$, where $g$ is given
   in Equation~\eqref{eq:gdefi}.
\end{cor}
\begin{rem} \label{rem:bbamend}
   Note that the statement of Corollary~\ref{cor:gnz} is wrong in the case $N_x = \{x\}$ for all $x \in \mcx$ if 
   $f$ is not hereditary and $\PP[\xi \neq 0] > 0$. As a counterexample consider the process that scatters a 
   fixed number $n \geq 1$ of points uniformly over $\mcx$ (cf.\ \cite{dvj88}, Example~14.2(a)). This process
   has a density given by $f(\varrho) = 
   e^{\alpha(\mcx)} \frac{n!}{\alpha(\mcx)^n} {\rm I}[\abs{\varrho} = n]$ and hence satisfies $g(x, \varrho) = 
   0$ for all $x \in \mcx$ and $\varrho \in \mfn$, which makes the right hand side in
   Corollary~\ref{cor:gnz} zero for every function $h$, whereas, with $h(x,\varrho) \equiv 1$, the
   left hand side is equal to $\EE \xi(\mcx) > 0$.\\
   Since the corollary does not hold generally, its use in the proofs of Theorem~2.4
   and~3.6 of \cite{bb92} and in the proof of Theorem~2.3 of \cite{chenxia04} is not justified unless an
   additional condition (such as hereditarity) is imposed.
\end{rem}
\begin{proof}[Proof of Proposition~\ref{prop:pregnz}]
   We proof the statement for non-negative $h$; the statement for bounded $h$ follows in the 
   usual way by decomposing $h$ into its positive and negative parts. The Slivnyak-Mecke theorem (see
   \cite{moellerwaage04}, Theorem~3.2, for $\mcx \subset \RD$, which clearly can be generalized to arbitrary 
   compact metric spaces) states that, for $\eta \sim P_1$ and measurable $\tilh: \mcx \times \mfn \to 
   \Rplus$,
   \begin{equation} \label{eq:sm}
      \EE \biggl( \int_{\mcx} \tilh(x, \eta - \delta_x) \; \eta(dx) \biggr) = \int_{\mcx} \EE \tilh(x, 
      \eta) \; \alpha(dx).
   \end{equation}
   Hence, setting $\tilh(x,\sigma) := h(x, \sigma \vert_{N_x^c}) f(\sigma+\delta_x)$, we obtain
   \begin{equation*}
   \begin{split}
      \int_{\mcx} \int_{\mfn(N_x^c)} h(x,\varrho) f_{N_x^c \cup \{x\}}(\varrho + \delta_x) \; 
         P_{1,N_x^c}(d \varrho) \; \alpha(dx) 
      &= \int_{\mcx} \int_{\mfn} h(x,\sigma \vert_{N_x^c}) f(\sigma +
         \delta_x) \; P_1(d \sigma) \; \alpha(dx) \\
      &= \int_{\mfn} \int_{\mcx} h(x,\sigma \vert_{N_x^c})
         f(\sigma) \; \sigma(dx) \; P_1(d\sigma) \\
      &= \EE \biggl( \int_{\mcx} h(x, \xi \vert_{N_x^c}) \; \xi(dx) \biggr).
   \end{split}
   \end{equation*}
   \vspace*{-7mm}

   \hspace*{2mm}
\end{proof}
\begin{proof}[Proof of Corollary~\ref{cor:gnz}] 
   The statement follows immediately from Proposition~\ref{prop:pregnz}, using that $f_{N_x^c \cup
   \{x\}}(\varrho + \delta_x) = g(x, \varrho) f_{N_x^c}(\varrho)$ for every $x \in \mcx$ and every $\varrho
   \in \mfn(N_x^c)$.
\end{proof}

\subsection{Density of the thinned process} \label{app:densityofthinned}

Let the point process $\xi$ and the random field $\pi$ be as for the definition of the thinning in Subsection~\ref{ssec:thinnings}. We assume additionally, as in Section~\ref{sec:bound}, that all the realizations of $\pi$ lie in an evaluable path space $E \subset [0,1]^{D}$ and that there is a regular conditional distribution of $\pi$ given the value of $\xi$ (see Appendix~\ref{app:lepsrcd}). It is essential for the construction below that we use the same such distribution throughout (i.e.\ without changing it in between on $\PP \xi^{-1}$-null sets), but insignificant, of course, which one we use.

Set then
\begin{equation*}
   q(\varrho \mvert \sigma) := \EE \biggl( \prod_{s \in \varrho} \pi(s) \prod_{\ts \in \sigma 
   \setminus \varrho} \bigl( 1 - \pi(\ts) \bigr) \biggm| \xi = \sigma \biggr)
\end{equation*}
for almost every $\sigma \in \mfn$ and for $\varrho \subset \sigma$.
It can be easily seen that the mapping $\bigl[ \mfn^2 \times E \to [0,1], (\varrho, \tvarrho, p) \mapsto \prod_{s \in \varrho} p(s) \prod_{\ts \in \tvarrho} \bigl( 1-p(\ts) \bigr) \bigr]$ is $\mcn^2 \otimes \mce$-measurable, whence we obtain that $q(\varrho \mvert \sigma)$ is well-defined and
that $\varphi: \mfn^2 \to [0,1], (\varrho, \tvarrho) \mapsto q(\varrho \mvert \varrho + \tvarrho \bigr)$ is measurable.

\begin{lem} \label{lem:fp}
   A density of the thinned process $\xi_{\pi}$ with respect to $P_{1}$ is given by
   \begin{equation*}
      f\obenpi(\varrho) := e^{\alpha(\mcx)} \int_{\mfn} q(\varrho \mvert \varrho + \tvarrho)
         f(\varrho + \tvarrho) \; P_{1}(d \tvarrho)
   \end{equation*}
   for almost every $\varrho \in \mfn$.
\end{lem}
\begin{proof}
The well-definedness and the measurability of $f\obenpi$ follow from the measurability of $\varphi$ defined above.

Consider two independent Poisson processes $\eta, \teta \sim P_{1}$. Take furthermore $\chi \sim P_{2}$ and let $\chi_{1/2}$ be a thinning of
$\chi$ with retention function $p \equiv 1/2$, which corresponds to picking a subset
of the points of $\chi$ uniformly at random. Note that \raisebox{0pt}[10.5pt][0pt]{$(\eta, \teta) \eqinlaw
(\chi_{1/2}, \chi \setminus \chi_{1/2})$} (see e.g.\ \cite{moellerwaage04}, Proposition~3.7, for $\mcx \subset \RD$; the proof can easily be adapted for general compact metric spaces).

Integration of the proposed density over an arbitrary set $C \in \mcn$, using Lemma~\ref{lem:hompoisdensity}
for the fifth line, yields
\begin{equation} \label{eq:showing_fp}
\begin{split}   
   \int_C f\obenpi(\varrho) \; P_{1}(d \varrho) &= e^{\alpha(\mcx)} \int_C
      \int_{\mfn} q(\varrho \mvert \varrho + \tvarrho) f(\varrho + \tvarrho)
      \; P_{1}(d \tvarrho) \; P_{1}(d \varrho) \\[2mm]
   &= e^{\alpha(\mcx)} \EE \Bigl( {\rm I} [ \eta \in C
      ] q(\eta \mvert \eta + \teta) f(\eta + \teta) \Bigr) \\[2.5mm]
   &= e^{\alpha(\mcx)} \EE \Bigl( {\rm I} [ \chi_{1/2} \in C ] q(\chi_{1/2} \mvert \chi) f
      (\chi) \Bigr) \\[2mm]
   &= e^{\alpha(\mcx)} \int_{\mfn} \frac{1}{2^{\abs{\sigma}}} \sum_{\varrho \subset \sigma}
      {\rm I} [ \varrho \in C ] q(\varrho \mvert \sigma) f(\sigma) \; P_{2}(d \sigma) \\
   &= \int_{\mfn} \sum_{\varrho \subset \sigma} {\rm I} [ \varrho \in C ] q(\varrho
      \mvert \sigma) f(\sigma) \; P_{1}(d \sigma) \\
   &= \int_{\mfn} \EE \bigl( Q_C\obenpi(\sigma) \bigm| \xi = \sigma \bigr) f(\sigma) \;
      P_1(d \sigma) \\
   &= \int_{\mfn} \PP [ \xi_{\pi} \in C \bigm| \xi =
      \sigma ] f(\sigma) \; P_{1}(d \sigma) \\[2mm]
   &= \PP [ \xi_{\pi} \in C ],
\end{split}
\end{equation}
where $Q_C\obenp(\sigma) := \sum_{\varrho \subset \sigma, \varrho \in C} \bigl( \prod_{s \in \varrho} p(s) \bigr) \bigl( \prod_{\ts \in \sigma \setminus \varrho} (1-p(\ts)) \bigr)$ for every $\sigma \in \mfn^{*}$ and every $p \in E$, so that $\bigl[(p, \sigma) \mapsto Q_C\obenp(\sigma)\bigr]$ is $\mce \otimes \mcn$-measurable and $Q_C\obenpi(\xi) = \PP[ \xi_{\pi}  \in C \mvert \xi, \pi ]$.
From Equation~\eqref{eq:showing_fp} the claim follows.
\end{proof}

The proof above yields that
\begin{equation*}
   \int_C f\obenpi(\varrho) \; P_{1}(d \varrho) = \int_{\mfn} \EE \bigl( Q_C\obenpi(\sigma) \mvert
   \xi = \sigma \bigr) f(\sigma) \; P_{1}(d \sigma)
\end{equation*}
for every $C \in \mcn$, and hence, by Equation~\eqref{eq:densityofrestricted}, that more generally, with $A \in \mcb$,
\begin{equation} \label{eq:wdt_crucial0}
\begin{split}
   \int_C f_{A} \obenpi(\varrho) \; P_{1,A}(d \varrho) &= \int_C
      \int_{\mfn(A^c)} f\obenpi
      (\varrho + \tvarrho) \; P_{1, A^c}(d \tvarrho) \; P_{1, A}(d \varrho) \\
   &= \int_{\tC} f\obenpi (\sigma) \; P_1(d \sigma) \\
   &= \int_{\mfn} \EE \bigl( Q_{\tC}\obenpi(\sigma) \bigm| \xi = 
      \sigma \bigr) f(\sigma) \; P_1(d \sigma) \\
   &= \int_{\mfn} \EE \bigl( Q_{C}\obenpi(\sigma \vert_{A})
      \bigm| \xi = \sigma \bigr) f(\sigma) \; P_1(d \sigma)
\end{split}
\end{equation}
for every $C \in \mcn(A)$, where $\tC := \{ \tsigma \in \mfn ; \, \tsigma \vert_A \in C \}$. For the last equality we used that $\sum_{\varrho \subset \sigma \vert_{A^c}} \bigl( \prod_{s \in \varrho} p(s) \bigr) \bigl( \prod_{\ts \in (\sigma \vert_{A^c}) \setminus \varrho} (1 - p(\ts)) \bigr) = 1$.

The various computations in \eqref{eq:showing_fp} remain correct (after the obvious minor modifications) if we add an extra point to $\varrho$, yielding in the unrestricted case
\begin{equation*}
\begin{split}
   \int_C f\obenpi(\varrho + \delta_x) \; P_{1}(d \varrho) &= \int_{\mfn} \sum_{\varrho \subset \sigma}
      {\rm I} [ \varrho \in C ] q(\varrho + \delta_x \mvert \sigma + \delta_x) f(\sigma + \delta_x) \; P_
      {1}(d \sigma) \\
   &= \int_{\mfn} \EE \bigl( \pi(x) Q_C\obenpi(\sigma) \bigm| \xi =
      \sigma + \delta_x \bigr) f(\sigma + \delta_x) \; P_{1}(d \sigma)
\end{split}
\end{equation*}
for every $C \in \mcn$, which holds for $\alpha$-almost every $x \in \mcx$.
Hence we obtain in a very similar fashion as in Equation~\eqref{eq:wdt_crucial0} that 
\begin{equation} \label{eq:wdt_crucial1}
   \int_C f_{A \cup \{x\}} \obenpi(\varrho + \delta_x) \; P_{1,A}(d \varrho) = 
      \int_{\mfn} \EE \bigl( \pi(x) Q_{C}\obenpi(\sigma \vert_{A})
      \bigm| \xi = \sigma + \delta_x \bigr) f(\sigma + \delta_x) \; P_1(d \sigma)
\end{equation}
for every $A \in \mcb$ and every $C \in \mcn(A)$, which holds for $\alpha$-almost every $x \in \mcx$.

\subsection{Technical conditions on $\bs{\xi}$ and $\bs{\pi}$: evaluable path space and regular conditional distribution}  \label{app:lepsrcd}

Consider a locally compact, second countable Hausdorff space $\mcy$ that is equipped with its Borel $\sigma$-algebra $\mcb = \mcb(\mcy)$. This is the most common type of space on which general point processes are defined. Any such space is separable, and a complete metric $\td$ can be introduced that generates its topology. With regard to the main part of this article, $\mcy$ is usually just our compact state space $\mcx$, but it is sometimes useful to consider a natural superset of $\mcx$ (e.g.\ $\RD$ if $\mcx \subset \RD$).
For sets of functions $\mcy \to [0,1]$, we introduce the concept of (locally) evaluable path spaces.
\begin{defi}
Let $E \subset [0,1]^{\mcy}$ and let $\mce$ be the canonical $\sigma$-algebra on $E$, which is generated by the evaluation mappings $\Psi_x: E \to [0,1], \, p \mapsto p(x)$, where $x \in \mcy$. For $U \in \mcb$
set furthermore $E(U) :=  \{ p \vert_U ; \, p \in E \}$ and write $\mce(U)$ for the corresponding $\sigma$-algebra generated by $\Psi_{U,x}: E(U) \to [0,1], \, \tp \mapsto \tp(x)$, where $x \in U$.
   \begin{enumerate}
   \item We call $E$ an \emph{evaluable path space} if the mapping $\Phi: E \times
      \mcy \to [0,1], (p,x) \mapsto p(x)$ is $\mce \otimes \mcb$-measurable.
   \item We call $E$ a \emph{locally evaluable path space} if the mapping
      $\Phi_U: E(U) \times U \to [0,1], (p,x) \mapsto p(x)$ is $\mce(U) \otimes \mcb_{U}$-measurable for
      every $U \subset \mcy$ that is open and relatively compact.
   \end{enumerate}
It can be easily seen that every locally evaluable path space is evaluable.
\end{defi}

For the main results of this article we assume that $\pi$ takes values in an evaluable path space
and that there exists a regular conditional distribution of $\pi$ given the value of $\xi$.
Neither of these assumptions presents a serious restriction, because they are both naturally satisfied in many practical applications, and if they are not, we can modify $\pi$ accordingly (provided it is measurable in the sense of Subsection~\ref{ssec:thinnings}) without changing the distribution of the resulting thinning. To see this let $R$ be a third of the minimal interpoint distance in $\xi$, which is positive except on a null set, and let $\tpi(\omega,x) := \pi(\omega,S(\omega))$ if there is a point $S(\omega)$ of $\xi(\omega)$ within distance $R(\omega)$ of $x$ and $\tpi(\omega,x) := 0$ otherwise. We have as path space $E$ for $\pi$ the space of all functions $p: \mcx \to [0,1]$ that are zero except on finitely many non-overlapping closed balls of positive radius, on each of which they are constant. By Proposition~\ref{prop:leps}(iii) below it can be seen, using the separability of $\mcy$, that this is an evaluable path space.
A regular conditional distribution of $\tpi$ given the value of $\xi$ can then be defined in a very straightforward manner, using the regular conditional distribution of $(\pi(s))_{s \in \sigma}$ given $\xi=\sigma$.

Since the above construction looks rather artificial in many situations, we provide a few manageable conditions under each of which a path space is (locally) evaluable, and hereby substantiate the statement that an evaluable path space is naturally present in many practical applications.
The proposition below is essentially the ``path space version'' of Proposition~A.D in \cite{schumi2}. Where it was conveniently possible, we have generalized the conditions from $\RplusD$ to the space $\mcy$.
\begin{defi}
   We call a set $E \subset [0,1]^{\mcy}$ \emph{separable from above [or below]} if there exists a countable
   set $\Sigma \subset \mcy$ such that for every $p \in E$, every open $\td$-ball $B \subset \mcy$ and every $y 
   \in \RR$ we have that $p(x) > y$ for all $x \in B \cap \Sigma$ implies $p(x) > y$ for all $x \in B$ [or
   $p(x) < y$ for all $x \in B \cap \Sigma$ implies $p(x) < y$ for all $x \in B$, respectively].
\end{defi}
\begin{prop}  \label{prop:leps}
   A set $E \subset [0,1]^{\mcy}$ is a locally evaluable path space if it satisfies any one of the 
   following conditions.
   \begin{enumerate}
      \item $\mcy = \RD$ and there is a closed convex cone $A \subset \RD$ of positive volume such
         that every $p \in E$ is continuous from $A$ (see the definition in \cite{schumi2}, Appendix~A.3);
      \item Every $p \in E$ is lower semicontinuous, and $E$ is separable from above;
      \item Every $p \in E$ is upper semicontinuous, and $E$ is separable from below;
      \item Every $p \in E$ is the indicator of a closed subset of $\mcy$, and the family $\mcc$ of these 
         closed subsets is separable in the sense that the subsets are jointly separable, i.e.\ there is a 
         countable set $\Sigma \subset \mcy$ such that $C = \overline{C \cap \Sigma}$ for every $C \in \mcc$,
         where the bar denotes topological closure.
   \end{enumerate}
\end{prop}
\begin{proof}
   (i) This follows with some minor adaptations from Proposition~A.D in $\cite{schumi2}$ by letting
   $\Omega := E$ 
   and defining $\pi$ to be the identity mapping on $\Omega$, so that $\pi(\omega,x) = \omega(x)$
   for every $\omega \in \Omega$ and every $x \in \RD$. Proposition~A.D(i) in \cite{schumi2} yields then that
   $\Phi_R: E(R) \times R \to [0,1], (p,x) \mapsto p(x)$ is $\mce(R) \otimes \mcb_R$-measurable
   for every bounded open rectangle $R \subset \RD$ (it is easy to see that the result from \cite{schumi2} 
   carries over from $\RplusD$ to $\RD$). The same result for a general bounded open set $U$ instead of $R$ 
   follows by writing $U$ as a countable union of open rectangles.

   (ii) We can essentially reproduce the prove of Proposition~A.D(ii) in $\cite{schumi2}$. Let $y \in
   [0,1)$ and 
   $U \subset \mcy$ be open and relatively compact. We show that $\Phi_U^{-1}\bigl((y,\infty)\bigr) \in
   \mce(U) \otimes \mcb_U$. Choose a separant $\Sigma$ as in the definition of the separability from above and 
   set $\mcg := \{ \BB^{\text{\raisebox{2pt}{$\ssst \hspace*{-0.7pt} \circ \!$}}}(x,1/n) \subset U; x \in 
   \Sigmap, n \in \NN \}$, where $\BB^{\text{\raisebox{2pt}{$\ssst \hspace*{-0.7pt} \circ \!$}}}(x,1/n)$ 
   denotes the open $\td$-ball with center $x$ and radius $1/n$, and $\Sigmap$ is an arbitrary countable dense 
   subset of $\mcy$.
   Noting that $p^{-1}\bigl((y,\infty)\bigr) \cap U$ is open by the lower semicontinuity of $p \in E$, we have
   \begin{equation*}
      \Phi_{U}^{-1}\bigl((y,\infty)\bigr) \, = \, \bigcup_{B \in \mcg} \biggl( \bigcap_{x \in B \cap \Sigma} 
      \Psi_{U,x}^{-1}\bigl((y,\infty)\bigr) \biggr) \times B \, \in \, \mce(U) \otimes \mcb_U.
   \end{equation*}

   (iii) Apply (ii) to $E' := \{ 1-p \, ; \: p \in E \}$.
   
   (iv) We apply (iii). It is evident that every $p \in E$ is upper semicontinuous. Separability of $E$ from
   below is inferred from the separability of $\mcc$ as follows. First note that the definition has to 
   be checked only for $y=1$, because $p \in E$ takes only the values $0$ and $1$. Let $\Sigma$ be a separant 
   for $\mcc$ as in statement (iv), and let $B$ be an open ball in $\mcy$. Then, for $C \in \mcc$ and
   $p = \indi{C}$, $p(x) < 1$ for all $x \in B \cap \Sigma$
   implies that $B \cap \Sigma \cap C = \emptyset$, hence $B \subset (C \cap \Sigma)^c$. Since
   $B$ is open, this implies $B \subset \interior \bigl( (C \cap \Sigma)^c \bigr) = \bigl( \overline{C
   \cap \Sigma} \bigr)^c = C^c$, where $\interior(A)$ denotes the interior of the set $A$ for any
   $A \subset \mcy$. Thus $p(x) < 1$ for all $x \in B$.
\end{proof}

\section*{Acknowledgement}

I would like to thank Adrian Baddeley for helpful discussions which have led to a considerable improvement of parts of this article. My further thanks go to the School of Mathematics and Statistics at the University of Western Australia for the hospitality I have received during my stay.

\bibliographystyle{plain}
\bibliography{thesis_db}

\end{document}